\documentclass[a4paper]{article}
\usepackage{amsfonts}
\usepackage{amsmath}
\usepackage{amssymb}
\begin{document}
\newtheorem{theorem}{Theorem}[section]
\newtheorem{lemma}[theorem]{Lemma}
\newtheorem{corollary}[theorem]{Corollary}
\newtheorem{conjecture}[theorem]{Conjecture}
\newtheorem{remark}[theorem]{Remark}
\newtheorem{definition}[theorem]{Definition}
\newtheorem{problem}[theorem]{Problem}
\newtheorem{example}[theorem]{Example}
\newtheorem{proposition}[theorem]{Proposition}
\title{{\bf PLURICANONICAL SYSTEMS OF PROJECTIVE VARIETIES OF GENERAL TYPE II}}
\date{November 29, 2004}
\author{Hajime TSUJI}
\maketitle
\begin{abstract}
We prove that there exists a positive integer $\nu_{n}$ depending
only on $n$ such that for every smooth projective $n$-fold of 
general type $X$ defined over complex numbers, 
$\mid mK_{X}\mid$ gives a birational rational map from $X$
into a projective space for every $m\geq \nu_{n}$. 
This theorem gives an affirmative answer to Severi's 
conjecture. MSC: 14J40, 32J18.
\end{abstract}
\tableofcontents
\section{Introduction}

Let $X$ be a smooth projective variety and let $K_{X}$ be the canonical 
bundle of $X$.
$X$ is said to be a general type, if there exists a positive 
integer $m$ such that the pluricanonical system 
$\mid mK_{X}\mid$ gives a birational (rational) embedding of $X$. 
The following problem is fundamental to study projective 
varieties of general type. \vspace{10mm}\\
{\bf Problem}
Find a positive integer $\nu_{n}$ depending only on $n$ 
such that for every smooth projective $n$-fold $X$,  
$\mid mK_{X}\mid$ gives a birational rational map from $X$
into a projective space for every $m \geqq \nu_{n}$. $\square$ \vspace{10mm} \\
If  $X$ is  a projective curve  of genus $\geqq 2$, it is well known that $\mid 3K_{X}\mid$ gives a 
projective embedding.
In the case that  $X$ is  a smooth projective surface  of general type, 
E. Bombieri showed that $\mid 5K_{X}\mid$ gives a birational
rational map from $X$ into a projective space (\cite{b3}).
But for the case of $\dim X\geqq 3$, very little is known about
the above problem.

The main purpose of this article is to prove the following theorems 
in full generality. 

\begin{theorem} 
There exists a positive integer $\nu_{n}$ which depends
only on $n$ such that for every smooth projective $n$-fold $X$
of general type defined over complex numbers, $\mid mK_{X}\mid$ gives a birational rational map
from $X$ into a projective space for every $m\geqq \nu_{n}$. $\square$
\end{theorem}

Theorem 1.1 is very much related to the theory of minimal models.
It has been conjectured that for every nonuniruled smooth projective
variety $X$, there exists a projective variety $X_{min}$ such that 
\begin{enumerate}
\item $X_{min}$ is birationally equivalent to $X$,
\item $X_{min}$ has only {\bf Q}-factorial terminal singularities,
\item $K_{X_{min}}$ is a nef {\bf Q}-Cartier divisor.
\end{enumerate}
$X_{min}$ is called a minimal model of $X$. 
To construct a minimal model, the minimal model program 
(MMP) has been proposed (cf. \cite[p.96]{k-m}).
The minimal model program was completed in the case of 
3-folds by S. Mori (\cite{mo}).

The proof of Theorem 1.1 can be very much simplified,
if we assume the existence of minimal models 
for projective varieties of general type (\cite{tu8}).
The proof here is modeled after the proof 
under the existence of minimal models  by using the theory of AZD
originated by the author (\cite{tu,tu2}). 

The major difficulty of the proof of 
Theorem 1.1 is to find {\bf ``a (universal) lower bound'' of the positivity 
of} $K_{X}$.  
In fact Theorem 1.1 is equivalent to the following theorem. 
\begin{theorem}
There exists a positive number $C_{n}$ which depends
only on $n$ such that for every smooth projective $n$-fold $X$
of general type defined over complex numbers, 
\[
\mu (X,K_{X}) := n!\cdot\overline{\lim}_{m\rightarrow\infty}m^{-\dim X}\dim H^{0}(X,{\cal O}_{X}(mK_{X})) \geqq C_{n}
\]
holds.
$\square$ \end{theorem} 
We note that $\mu (X,K_{X})$ is equal to the intersection 
number $K_{X}^{n}$ for a minimal projective $n$-fold
$X$ of general type. 
In Theorems 1.1 and 1.2, the numbers $\nu_{n}$ and 
$C_{n}$ have not yet been computed effectively. 

The relation between Theorems 1.1 and 1.2 is as follows.
Theorem 1.2 means that there exists a universal lower bound of 
the positivity of canonical bundle of smooth projective variety of 
general type with a fixed dimension.
On the other hand, for a smooth projective variety of general type $X$, 
 let us consider the lower bound of $m$ such that $\mid mK_{X}\mid$ gives a birational embedding.
Such a lower bound depends on the positivity of $K_{X}$ on certain 
 subvarieties which appear as the strata of the filtrations as in 
\cite{t,a-s}(cf. Section 3.1).

The positivity of $K_{X}$ on the subvarieites  can be related 
to the positivity of the canonical bundles of the smooth models 
of the subvarieties via the subadjunction theorem due to Kawamata (\cite{ka}). 
We note that for a smooth projective variety $X$ of general type
there exists a nonempty open subset $U_{0}$ in countable Zariski topology such that for every $x\in U_{0}$, any subvariety containing $x$ is of general type. 

The organization of the paper is as follows. 
In Section 2, we review the basic techniques to prove Theorems 1.1 and 1.2. 

In Section 3, we prove Theorems 1.1 and 1.2 without assuming the 
existence of minimal models for projective varieties of general type.
Here we use the AZD (cf. Section 2.2) of $K_{X}$ instead of minimal models. 
And we use the subadjunction theorem (Theorem \ref{subad1}) 
and the positivity theorem (Theorem \ref{pos}) due to Kawamata. 

In Section 4, we discuss the application of Theorems 1.1 and 1.2 
to Severi-Iitaka's conjecture. 

In this paper all the varieties are defined over {\bf C}.

This is the continuation of the paper \cite{tu8} and is a transcription 
of the latter half of \cite{tu7}.

The author would like to express his sincere thanks to Professor Akira Fujiki 
who helped him to improve the exposition.

\section{Preliminaries}

In this section, we shall summerize the basic analyic tools to  
prove Theorems 1.1 and 1.2 by transcripting the proof of Theorems 1.1 and 1.2
assuming MMP (\cite{tu8}).
   
\subsection{Multiplier ideal sheaves and singularities of divisors}

In this subsection, we shall review the 
relation between multiplier ideal sheaves and singularities of divisors.
Throughout this subsection $L$ will denote a holomorphic line bundle on a complex manifold $M$. 
\begin{definition}
A  singular hermitian metric $h$ on $L$ is given by
\[
h = e^{-\varphi}\cdot h_{0},
\]
where $h_{0}$ is a $C^{\infty}$-hermitian metric on $L$ and 
$\varphi\in L^{1}_{loc}(M)$ is an arbitrary function on $M$.
We call $\varphi$ the  weight function of $h$ with respect to $h_{0}$.
$\square$ \end{definition}
The curvature current $\Theta_{h}$ of the singular hermitian line
bundle $(L,h)$ is defined by
\[
\Theta_{h} := \Theta_{h_{0}} + \sqrt{-1}\partial\bar{\partial}\varphi ,
\]
where $\partial\bar{\partial}$ is taken in the sense of a current.
The $L^{2}$-sheaf ${\cal L}^{2}(L,h)$ of the singular hermitian
line bundle $(L,h)$ is defined by
\[
{\cal L}^{2}(L,h)(U) := \{ \sigma\in\Gamma (U,{\cal O}_{M}(L))\mid 
\, h(\sigma ,\sigma )\in L^{1}_{loc}(U)\} ,
\]
where $U$ runs over the  open subsets of $M$.
In this case there exists an ideal sheaf ${\cal I}(h)$ such that
\[
{\cal L}^{2}(L,h) = {\cal O}_{M}(L)\otimes {\cal I}(h)
\]
holds.  We call ${\cal I}(h)$ the {\bf multiplier ideal sheaf} of $(L,h)$.
If we write $h$ as 
\[
h = e^{-\varphi}\cdot h_{0},
\]
where $h_{0}$ is a $C^{\infty}$ hermitian metric on $L$ and 
$\varphi\in L^{1}_{loc}(M)$ is the weight function, we see that
\[
{\cal I}(h) = {\cal L}^{2}({\cal O}_{M},e^{-\varphi})
\]
holds.
For $\varphi\in L^{1}_{loc}(M)$ we define the multiplier ideal sheaf of $\varphi$ by 
\[
{\cal I}(\varphi ) := {\cal L}^{2}({\cal O}_{M},e^{-\varphi}).
\]
\begin{example}\label{ex}
Let $m$ be a positive integer. 
Let $\sigma\in \Gamma (X,{\cal O}_{X}(mL))$ be the global section. 
Then 
\[
h := \frac{1}{\mid\sigma\mid^{2/m}} = \frac{h_{0}}{(h_{0}^{m}(\sigma ,\sigma))^{1/m}}
\]
is a singular hemitian metric on $L$, 
where $h_{0}$ is an arbitrary $C^{\infty}$-hermitian metric on $L$
(the righthand side is ovbiously independent of $h_{0}$).
The curvature $\Theta_{h}$ is given by
\[
\Theta_{h} = \frac{2\pi\sqrt{-1}}{m}(\sigma ),
\]
where $(\sigma )$ denotes the current of integration over the 
divisor of $\sigma$. $\square$ 
\end{example}
\begin{definition}
$L$ is said to be pseudoeffective, if there exists 
a singular hermitian metric $h$ on $L$ such that 
the curvature current 
$\Theta_{h}$ is a closed positive current.

Also a singular hermitian line bundle $(L,h)$ is said to be pseudoeffective, 
if the curvature current $\Theta_{h}$ is a closed positive current.
$\square$ \end{definition}

Let $m$ be a positive integer and  $\{\sigma_{i}\}$ a finite number of global holomorphic sections of $mL$. 
Let  $\phi$ be a $C^{\infty}$-function on $M$.
Then 
\[
h := e^{-\phi}\cdot\frac{1}{\sum_{i}\mid\sigma_{i}\mid^{2/m}}
\]
defines a singular hermitian metric  on 
$L$.
We call such a metric $h$ a singular hermitian metric 
on $L$ with  {\bf algebraic singularities}.
Singular hermitian metrics with algebraic singularities 
are particulary easy to handle, because its multiplier 
ideal sheaf of the metric can 
be controlled by taking a resolution  
of the base scheme  $\cap_{i}(\sigma_{i})$.

Let $D= \sum a_{i}D_{i}$ be an effective {\bf Q}-divisor on $X$. 
Let $\sigma_{i}$ be a section of ${\cal O}_{X}(D_{i})$ with divisor $D_{i}$
respectively. 
Then we define 
\[
{\cal I}(D) : = {\cal I}(\sum_{i}a_{i}\log h_{i}(\sigma_{i},\sigma_{i}))
\]
and call it the multiplier ideal sheaf of the divisor $D$, where 
$h_{i}$ denotes a $C^{\infty}$-hermitian metric 
of ${\cal O}_{X}(D_{i})$ respectively.   
It is clear that ${\cal I}(D)$ is independent of the choice of
the hermitian metrics $\{ h_{i}\}$.

Let us consider the relation between ${\cal I}(D)$ and 
singularities of $D$. 
\begin{definition}\label{KLT}
Let $X$ be a normal variety and $D= \sum_{i}d_{i}D_{i}$ an effective {\bf Q}-divisor such that $K_{X}+D$ is {\bf Q}-Cartier. 
If $\mu : Y \longrightarrow X$ is a log resolution  of the pair 
$(X,D)$, i.e., $\mu$ is a composition of successive blowing ups with smooth centers 
such that $Y$ is smooth and the support of $f^{*}D$ is a divisor with normal crossings, then we can write
\[
K_{Y} + \mu_{*}^{-1}D = \mu^{*}(K_{X}+D) + F
\]
with $F = \sum_{j}e_{j}E_{j}$ for the exceptional divisors $\{ E_{j}\}$, 
where $\mu_{*}^{-1}D$ denotes the strict transform of $D$. 
We call $F$ the discrepancy and $e_{j}\in  \mbox{\bf Q}$ the discrepancy
coefficient for $E_{j}$. 
We regard $-d_{i}$ as the discrepancy coefficient of $D_{i}$. 

The pair $(X,D)$ is said to have only {\bf Kawamata log terminal singularities}
({\bf KLT})(resp. {\bf log canonical singularities}({\bf LC})), if 
$d_{i}< 1$(resp. $\leqq 1$) for all $i$ and $e_{j} > -1$ (resp. $\geqq -1$)
for all $j$ for a log resolution $\mu : Y \longrightarrow X$.
One can also say that $(X,D)$ is KLT (resp. LC), or $K_{X}+D$ is KLT
(resp. LC), when $(X,D)$ has only KLT (resp. LC).
The pair $(X,D)$ is said to be KLT (resp. LC) at a point $x_{0}\in X$,
if $(U,D\mid_{U})$ is KLT (resp. LC) for some neighbourhood $U$ 
of $x_{0}$. 
$\square$ \end{definition}
The following proposition is a dictionary between algebraic geometry and 
the $L^{2}$-method.

\begin{proposition}\label{prop}
Let $D$ be a divisor on a smooth projective variety $X$. 
Then $(X,D)$ is KLT, if and only if 
${\cal I}(D)$ is trivial ($= {\cal O}_{X}$).
$\square$ \end{proposition}
The proof is trivial and left to the reader. 
To locate the co-support of the multiplier ideal the following notion
is useful.  
\begin{definition}\label{CLC}
A subvariety $W$ of $X$ is said to be a {\bf center of log canonical singularities} for the pair $(X,D)$, if there is a birational morphism from 
a normal variety $\mu : Y \longrightarrow X$ and a prime divisor $E$ on $Y$
with the discrepancy coefficient $e \leqq -1$
such that $\mu (E) = W$. 
$\square$ \end{definition}
The set of all the centers of log canonical singularities is denoted 
by $CLC(X,D)$.
For a point $x_{0}\in X$, we define
$CLC(X,x_{0},D) := \{ W\in CLC(X,D)\mid x_{0}\in W\}$.
We quote the following proposition to introduce the notion of 
the minimal center of logcanoical singularities. 

\begin{proposition}\label{MC}(\cite[p.494, Proposition 1.5]{ka2})
Let $X$ be a  normal variety and $D$ an effective {\bf Q}-Cartier divisor 
such that $K_{X}+D$ is {\bf Q}-Cartier. 
Assume that $X$ is KLT and $(X,D)$ is LC.
If $W_{1},W_{2}\in CLC(X,D)$ and $W$ an irreducible component of $W_{1}\cap W_{2}$, then $W\in CLC(X,D)$. 
This implies that if $(X,D)$ is LC but not KLT, 
then there exists the unique minimal element of $CLC(X,D)$. 
Also if $(X,D)$ is LC but not KLT at a point $x_{0}\in X$,
then there exists the unique minimal element of $CLC(X,x_{0},D)$. 
\end{proposition}
We call these  minimal elements  the {\bf minimal center of LC singularities}
of $(X,D)$ and the {\bf minimal center of LC singularities  
of $(X,D)$ at $x_{0}$} respectively.

\subsection{Analytic Zariski decomposition}

To study a pseudoeffective line bundle we introduce the notion of analytic Zariski
decompositions.
By using analytic Zariski decompositions, we can handle a pseudoeffective line bundle, as if it  were a nef line bundle.
\begin{definition}\label{AZDdef}
Let $M$ be a compact complex manifold and let $L$ be a line bundle
on $M$.  A singular hermitian metric $h$ on $L$ is said to be 
an {\bf analytic Zariski decomposition} ({\bf AZD} in short), if the followings hold.
\begin{enumerate}
\item $\Theta_{h}$ is a closed positive current,
\item for every $m\geqq 0$, the natural inclusion
\[
H^{0}(M,{\cal O}_{M}(mL)\otimes{\cal I}(h^{m}))\rightarrow
H^{0}(M,{\cal O}_{M}(mL))
\]
is isomorphim.
\end{enumerate}
$\square$ \end{definition}
\begin{remark} If an AZD exists on a line bundle $L$ on a smooth projective
variety $M$, $L$ is pseudoeffective by the condition 1 above.
$\square$ \end{remark}

\begin{theorem}\label{AZD1}(\cite{tu,tu2})
 Let $L$ be a big line  bundle on a smooth projective variety
$M$.  Then $L$ has an AZD. 
$\square$ \end{theorem}
As for the existence for general pseudoeffective line bundles, 
now we have the following theorem.
\begin{theorem}\label{AZD}(\cite[Theorem 1.5]{d-p-s})
Let $X$ be a smooth projective variety and let $L$ be a pseudoeffective 
line bundle on $X$.  Then $L$ has an AZD.
$\square$ \end{theorem}
Although the proof is in \cite{d-p-s}, 
we shall give a proof here, because we shall use it afterward. 

 Let  $h_{0}$ be a fixed $C^{\infty}$-hermitian metric on $L$.
Let $E$ be the set of singular hermitian metric on $L$ defined by
\[
E = \{ h ; h : \mbox{lowersemicontinuous singular hermitian metric on $L$}, 
\]
\[
\hspace{70mm}\Theta_{h}\,
\mbox{is positive}, \frac{h}{h_{0}}\geq 1 \}.
\]
Since $L$ is pseudoeffective, $E$ is nonempty.
We set 
\[
h_{L} = h_{0}\cdot\inf_{h\in E}\frac{h}{h_{0}},
\]
where the infimum is taken pointwise. 
The supremum of a family of plurisubharmonic functions 
uniformly bounded from above is known to be again plurisubharmonic, 
if we modify the supremum on a set of measure $0$(i.e., if we take the uppersemicontinuous envelope) by the following theorem of P. Lelong.

\begin{theorem}\label{Lelong}(\cite[p.26, Theorem 5]{l})
Let $\{\varphi_{t}\}_{t\in T}$ be a family of plurisubharmonic functions  
on a domain $\Omega$ 
which is uniformly bounded from above on every compact subset of $\Omega$.
Then $\psi = \sup_{t\in T}\varphi_{t}$ has a minimum 
uppersemicontinuous majorant $\psi^{*}$  which is plurisubharmonic.
We call $\psi^{*}$ the uppersemicontinuous envelope of $\psi$. 
$\square$ \end{theorem}
\begin{remark} In the above theorem the equality 
$\psi = \psi^{*}$ holds outside of a set of measure $0$(cf.\cite[p.29]{l}). 
$\square$ \end{remark}

By Theorem \ref{Lelong}, we see that $h_{L}$ is also a 
singular hermitian metric on $L$ with $\Theta_{h_{L}}\geq 0$.
Suppose that there exists a nontrivial section 
$\sigma\in \Gamma (X,{\cal O}_{X}(mL))$ for some $m$ (otherwise the 
second condition in Definition \ref{AZDdef} is empty).
We note that  
\[
\frac{1}{\mid\sigma\mid^{\frac{2}{m}}} 
\]
gives the weight of a singular hermitian metric on $L$ with curvature 
$2\pi m^{-1}(\sigma )$, where $(\sigma )$ is the current of integration
along the zero set of $\sigma$. 
By the construction we see that there exists a positive constant 
$c$ such that  
\[
\mbox{($*$)}\hspace{10mm} \frac{h_{0}}{\mid\sigma\mid^{\frac{2}{m}}} \geq c\cdot h_{L}
\]
holds. 
Hence
\[
\sigma \in H^{0}(X,{\cal O}_{X}(mL)\otimes{\cal I}_{\infty}(h_{L}^{m}))
\]
holds.  
In praticular
\[
\sigma \in H^{0}(X,{\cal O}_{X}(mL)\otimes{\cal I}(h_{L}^{m}))
\]
holds.  
 This means that $h_{L}$ is an AZD of $L$. 
\vspace{10mm} $\square$ 

\begin{remark}\label{r2.3}
By the above proof (see ($*$)) we have that for the AZD $h_{L}$ constructed 
as above
\[
H^{0}(X,{\cal O}_{X}(mL)\otimes{\cal I}_{\infty}(h_{L}^{m}))
\simeq 
H^{0}(X,{\cal O}_{X}(mL))
\]
holds for every $m$, where ${\cal I}_{\infty}(h_{L}^{m})$ 
denotes the $L^{\infty}$-multiplier ideal sheaf, i.e., 
for every open subset $U$ in $X$, 
\[
{\cal I}_{\infty}(h_{L}^{m})(U):= \{ f\in {\cal O}_{X}(U) \mid
\mid f\mid^{2}(h_{L}/h_{0})^{m}\in L^{\infty}_{\mbox{loc}}(U)\} ,
\]
where $h_{0}$ is a $C^{\infty}$-hermitian metric on $L$. 
$\square$ \end{remark}

Entirely the same proof as that of Theorem \ref{AZD}, we obtain the following 
corollary. 
\begin{corollary}\label{cor}
Let $(L,h_{0})$ be a singular hermitian line bundle on a compact K\"{a}hler  
manifold $(X,\omega )$.  Suppose that 
\[
E(L,h_{0}):= \{ \varphi\in L^{1}_{loc}(X) \mid  \varphi \leqq 0, \,\,\Theta_{h_{0}} + \sqrt{-1}\partial\bar{\partial}\varphi \geqq 0 \}
\]
is nonempty.   
Then if we define the function $\varphi_{P}\in L^{1}_{loc}(X)$ by
\[
\varphi_{P}(x)  := \sup \{\varphi (x)\mid \varphi \in E\} \,\,\,\, (x\in X).
\]
Then $h := e^{-\varphi_{P}}\cdot h_{0}$ is a singular hermitian metric on 
$L$ such that 
\begin{enumerate}
\item $\Theta_{h}\geqq 0$. 
\item $H^{0}(X,{\cal O}_{X}(mL)\otimes {\cal I}_{\infty}(h^{m})) 
\simeq H^{0}(X,{\cal O}_{X}(mL)\otimes {\cal I}_{\infty}(h_{0}^{m}))$ holds 
for every $m\geqq 0$. $\square$
\end{enumerate}
\end{corollary}
We call $h$ an AZD of $(L,h_{0})$.  This is a slight generalization of 
the notion of 
AZD's of  pseudoeffective line bundles. 

\begin{remark}\label{r2.4}
In Corollary \ref{cor}, $E(L,h_{0})$ is nonempty, if there exists a positive integer 
$m_{0}$ and $\sigma  \in  H^{0}(X,{\cal O}_{X}(m_{0}L)\otimes {\cal I}_{\infty}(h_{0}^{m_{0}}))$ such that $h_{0}^{m_{0}}(\sigma ,\sigma ) \leqq 1$. 
In this case 
\[
\varphi := \frac{1}{m_{0}}\log h_{0}^{m_{0}}(\sigma ,\sigma )
\]
belongs to $E(L,h_{0})$. 
\end{remark}

\subsection{The $L^{2}$-extension theorem}

 Let $M$ be a complex manifold of dimension $n$ and let $S$ be a closed complex submanifold of $M$. 
Then we consider a class of continuous function $\Psi : M\longrightarrow [-\infty , 0)$  such that  
\begin{enumerate}
\item $\Psi^{-1}(-\infty ) \supset S$,
\item if $S$ is $k$-dimensional around a point $x$, there exists a local 
coordinate system $(z_{1},\ldots ,z_{n})$ on a neighbourhood of $x$ such that 
$z_{k+1} = \cdots = z_{n} = 0$ on $S\cap U$ and 
\[
\sup_{U\backslash S}\mid \Psi (z)-(n-k)\log\sum_{j=k+1}^{n}\mid z_{j}\mid^{2}\mid < \infty .
\]
\end{enumerate} 
The set of such functions $\Psi$ will be denoted by $\sharp (S)$. 

For each $\Psi \in \sharp (S)$, one can associate a positive measure 
$dV_{M}[\Psi ]$ on $S$ as the minimum element of the 
partially ordered set of positive measures $d\mu$ 
satisfying 
\[
\int_{S_{k}}f\, d\mu \geqq 
\overline{\lim}_{t\rightarrow\infty}\frac{2(n-k)}{v_{2n-2k-1}}
\int_{M}f\cdot e^{-\Psi}\cdot \chi_{R(\Psi ,t)}dV_{M}
\]
for any nonnegative continuous function $f$ with 
$\mbox{supp}\, f\subset\subset M$.
Here $S_{k}$ denotes the $k$-dimensional component of $S$,
$v_{m}$ denotes the volume of the unit sphere 
in $\mbox{\bf R}^{m+1}$ and 
$\chi_{R(\Psi ,t)}$ denotes the characteristic funciton of the set 
\[
R(\Psi ,t) = \{ x\in M\mid -t-1 < \Psi (x) < -t\} .
\]

Let $M$ be a complex manifold and let $(E,h_{E})$ be a holomorphic hermitian vector 
bundle over $M$. 
Given a positive measure $d\mu_{M}$ on $M$,
we shall denote $A^{2}(M,E,h_{E},d\mu_{M})$ the space of 
$L^{2}$ holomorphic sections of $E$ over $M$ with respect to $h_{E}$ and 
$d\mu_{M}$. 
Let $S$ be a closed  complex submanifold of $M$ and let $d\mu_{S}$ 
be a positive measure on $S$. 
The measured submanifold $(S,d\mu_{S})$ is said to be a set of 
interpolation for $(E,h_{E},d\mu_{M})$, or for the 
sapce $A^{2}(M,E,h_{E},d\mu_{M})$, if there exists a bounded linear operator
\[
I : A^{2}(S,E\mid_{S},h_{E},d\mu_{S})\longrightarrow A^{2}(M,E,h_{E},d\mu_{M})
\]
such that $I(f)\mid_{S} = f$ for any $f\in A^{2}(S,E\mid_{S},h_{E},d\mu_{S})$. 
$I$ is called an interpolation operator.
The following theorem is crucial.

\begin{theorem}\label{l2}(\cite[Theorem 4]{o})
Let $M$ be a complex manifold with a continuous volume form $dV_{M}$,
let $E$ be a holomorphic vector bundle over $M$ with $C^{\infty}$-fiber 
metric $h_{E}$, let $S$ be a closed complex submanifold of $M$,
let $\Psi\in \sharp (S)$ and let $K_{M}$ be the canonical bundle of $M$.
Then $(S,dV_{M}(\Psi ))$ is a set of interpolation 
for $(E\otimes K_{M},h_{E}\otimes (dV_{M})^{-1},dV_{M})$, if 
the followings are satisfied.
\begin{enumerate}
\item There exists a closed set $X\subset M$ such that 
\begin{enumerate}
\item $X$ is locally negligble with respect to $L^{2}$-holomorphic functions, i.e., 
for any local coordinate neighbourhood $U\subset M$ and for any $L^{2}$-holomorphic function $f$ on $U\backslash X$, there exists a holomorphic function 
$\tilde{f}$ on $U$ such that $\tilde{f}\mid U\backslash X = f$.
\item $M\backslash X$ is a Stein manifold which intersects with every component of $S$. 
\end{enumerate}
\item $\Theta_{h_{E}}\geqq 0$ in the sense of Nakano,
\item $\Psi \in \sharp (S)\cap C^{\infty}(M\backslash S)$,
\item $e^{-(1+\epsilon )\Psi}\cdot h_{E}$ has semipositive 
curvature in the sense of Nakano for every $\epsilon \in [0,\delta]$ 
for some $\delta > 0$.
\end{enumerate}
Under these conditions, there exists a constant $C$ and an interpolation operator 
from $A^{2}(S,E\otimes K_{M}\mid_{S},h\otimes (dV_{M})^{-1}\mid_{S},dV_{M}[\Psi ])$
to $A^{2}(M,E\otimes K_{M},h\otimes (dV_{M})^{-1}.dV_{M})$ whose 
norm does not exceed $C\cdot\delta^{-3/2}$.
If $\Psi$ is plurisubharmonic, the interpolation operator can be chosen 
so that its norm is less than $2^{4}\pi^{1/2}$.
$\square$ \end{theorem}

The above theorem can be generalized to the case that 
$(E,h_{E})$ is a singular hermitian line bundle with semipositive
curvature current  (we call such a singular hermitian line 
bundle $(E,h_{E})$ a {\bf pseudoeffective singular hermitian line bundle}) as was remarked in \cite{o}. 

\begin{lemma}\label{l2lemma} 
Let $M,S,\Psi ,dV_{M}, dV_{M}[\Psi], (E,h_{E})$ be as in Theorem \ref{l2}. 
Let $(L,h_{L})$ be a pseudoeffective singular hermitian line 
bundle on $M$. 
Then $(S,dV_{M}[\Psi ])$ is a set of interpolation for 
$(K_{M}\otimes E\otimes L,dV_{M}^{-1}\otimes h_{E}\otimes h_{L})$.  
$\square$ \end{lemma}

\subsection{A construction of the function $\Psi$}
Here we shall show the standard construction of the function $\Psi$ in 
Theorem \ref{l2}.
Let $M$ be a smooth projective $n$-fold and 
let $S$ be a $k$-dimensional (not necessary smooth)
subvariety of $M$. 
Let ${\cal U} = \{ U_{\gamma}\}$ be a finite 
Stein covering of $M$ and 
let $\{ f^{(\gamma )}_{1},\ldots , f_{m(\gamma )}^{(\gamma )}\}$ 
be a generator of the ideal sheaf associated with $S$ 
on $U_{\gamma}$. 
Let $\{ \phi_{\gamma}\}$ be a partition of unity which  subordinates
to ${\cal U}$. 
We set 
\[
\Psi := (n-k)\sum_{\gamma}\phi_{\gamma}\cdot (\sum_{\ell = 1}^{m(\gamma )}
\mid f_{\ell}^{(\gamma )}\mid^{2}).
\]
Then the residue volume form $dV[\Psi ]$ is defined 
as in the last subsection. 
Here the residue volume form $dV[\Psi ]$  
of a continuous volume form $dV$ on $M$ is not well defined on the singular
locus of $S$.
But this is not a difficulty to apply Theorem \ref{l2} or 
Lemma \ref{l2lemma}, since there exists a proper Zariski closed subset  
$Y$ of $X$ such that 
$(X - Y)\cap S$ is smooth. 

\subsection{Volume of pseudoeffective line bundles}

To measure the positivity of big line bundles on a projective 
variety, we shall introduce the notion of volume of a projective 
variety with respect to a big line bundle. 

\begin{definition}\label{volume1} Let $L$ be a line bundle on a compact complex 
manifold $M$ of dimension $n$. 
We define the {\bf volume} of $M$ with respect to $L$ by
\[
\mu (M,L) := n!\cdot\overline{\lim}_{m\rightarrow\infty}m^{-n}
\dim H^{0}(M,{\cal O}_{M}(mL)).
\]
$\square$ \end{definition}
With respect to a pseudoeffective singular hermitian line bundle
(for the definition of  pseudoeffective singular hermitian line bundles, 
see the last part of Section 2.3), 
we define the volume as follows. 
\begin{definition}\label{volume}(\cite{tu3})
Let $(L,h)$ be a pseudoeffective singular hermitian line bundle on a smooth projective variety
$X$ of dimension $n$. 
We define the {\bf volume of $X$ with respect to 
$(L,h)$ }by 
\[
\mu (X,(L,h)) := n!\cdot\overline{\lim}_{m\rightarrow\infty}m^{-n}
\dim H^{0}(X,{\cal O}_{X}(mL)\otimes{\cal I}(h^{m})).
\]
A pseudoeffective singular hermitian line bundle $(L,h)$ 
is said to be big, if $\mu (X,(L,h)) > 0$ holds. 

We may consider $\mu (X,(L,h))$ as the {\bf intersection number}  
$(L,h)^{n}$. 
We also denote $\mu (X,(L,h))$ by $(L,h)^{n}$. 
Let  $Y$ be  a subvariety of $X$ of dimension $d$
and let $\pi_{Y} : \tilde{Y}\longrightarrow Y$ be 
a resolution of $Y$.   We define 
$\mu (Y,(L,h)\mid_{Y})$ as 
\[
\mu (Y,(L,h)\mid_{Y})
:= \mu (\tilde{Y},\pi^{*}_{Y}(L,h)).
\]
The righthand side is independent of the choice of the 
resolution $\pi$ because of the remark below. 
We also denote $\mu (Y,(L,h)\mid_{Y})$ by $(L,h)^{d}\cdot Y$.
$\square$ \end{definition}
\begin{remark}
Let us use the same notations in Definition \ref{volume}.
Let $\pi : \tilde{X}\longrightarrow X$ 
be any modification. 
Then 
\[
\mu (X,(L,h)) = \mu (\tilde{X},\pi^{*}(L,h))
\]
holds, since
\[
\pi_{*}({\cal O}_{\tilde{X}}(K_{\tilde{X}})\otimes {\cal I}(\pi^{*}h^{m}))
= {\cal O}_{X}(K_{X})\otimes {\cal I}(h^{m})
\]
holds for every $m$ and 
\[
\overline{\lim}_{m\rightarrow\infty}m^{-n}
\dim H^{0}(X,{\cal O}_{X}(mL)\otimes{\cal I}(h^{m}))
=  
\overline{\lim}_{m\rightarrow\infty}m^{-n}
\dim H^{0}(X,{\cal O}_{X}(mL + D)\otimes{\cal I}(h^{m}))
\]
holds for any Cartier divisor $D$ on $X$. 
This last equality can be easily checked, if $D$ is a smooth 
irreducible divisor, by using 
the exact sequence 
\[
0 \rightarrow {\cal O}_{X}(mL)\otimes {\cal I}(h^{m})
\rightarrow {\cal O}_{X}(mL+D)\otimes {\cal I}(h^{m})
\rightarrow {\cal O}_{D}(mL+D)\otimes {\cal I}(h^{m})
\rightarrow 0.
\]
For a general $D$, the equality follows by expressing $D$ 
as a difference of two very ample divisors. 
$\square$ \end{remark}

\subsection{A subadjunction theorem}
Let $M$ be a smooth projective variety 
and let $(L,h_{L})$ be a singular hermitian line bundle on $M$ such that 
$\Theta_{h_{L}}\geqq 0$ on $M$.  
We assume that $h_{L}$ is lowersemicontinuous. 
This is a technical assumption so that a local potential 
of the curvature current of $h$ is plurisubharmonic. 

Let $m_{0}$ be a positive integer. 
Let $\sigma \in \Gamma (M,{\cal O}_{M}(m_{0}L)\otimes {\cal I}(h))$ be a 
global section. 
Let $\alpha$ be a positive rational number $\leqq 1$ and let $S$ be 
an irreducible subvariety of $M$ 
such that  $(M, \alpha (\sigma ))$ is LC(log canonical) but not KLT(Kawamata log terminal)
on the generic point of $S$ and $(M,(\alpha -\epsilon )(\sigma ))$ is KLT on the generic point of $S$ 
for every $0 < \epsilon << 1$. 
We set 
\[
\Psi_{S} = \alpha \log h_{L}(\sigma ,\sigma ).
\]
Suppose that $S$ is smooth for simplicity 
(if $S$ is not smooth, we just need to take an embedded 
resolution to apply Theorems \ref{subad1}, \ref{subad2} below). 
We shall assume that $S$ is not contained in the 
singular locus of $h_{L}$, where the singular locus of $h_{L}$ means the 
set of points where $h$ is $+\infty$. 
Let $dV$ be a $C^{\infty}$-volume form on $M$. 

Then as in Section 2.3, we may define a (possibly singular) measure 
$dV[\Psi_{S} ]$ on $S$. 
This can be viewed as follows. 
Let $f : N \longrightarrow M$ be a log resolution of 
$(X,\alpha (\sigma ))$. 
Then  as in Section 2.4,
 we may define the singular volume form $f^{*}dV[f^{*}\Psi_{S} ]$ 
on the divisorial component of $f^{-1}(S)$ (the volume form is identically 
$0$ on the components with discrepancy $> -1$).
The singular volume form $dV[\Psi_{S} ]$ is defined as the fibre integral of 
$f^{*}dV[f^{*}\Psi_{S} ]$ (the actual integration takes place only on the components with discrepancy $-1$). 
Let $d\mu_{S}$ be a $C^{\infty}$-volume form on $S$ and 
let $\varphi$ be the function on $S$ defined by
\[
\varphi := \log \frac{d\mu_{S}}{dV[\Psi_{S} ]}
\]
($dV[\Psi_{S} ]$ may be singular on a subvariety of $S$, also 
it may be totally singular on $S$). 

\begin{theorem}\label{subad1}(\cite[Theorem 5.1]{tu5})
Let $M$,$S$,$\Psi_{S}$ be as above. 
Suppose that $S$ is smooth.   
Let $d$ be a positive integer such that $d > \alpha m_{0}$. 
Then every element of 
$A^{2}(S,{\cal O}_{S}(m(K_{M}+dL)),e^{-(m-1)\varphi}\cdot dV^{-m}\cdot h_{L}^{m}\mid_{S},dV[\Psi_{S} ])$ 
extends to an element of 
\[
H^{0}(M,{\cal O}_{M}(m(K_{M}+dL))). 
\]
$\square$ \end{theorem}
As we mentioned as above the smoothness assumption on $S$ is 
just to make the statement simpler.  

Theorem \ref{subad1} follows from 
Theorem \ref{subad2} below by  minor modifications (cf. \cite{tu5}). 
The main difference is the fact that the residue volume form 
$dV[\Psi_{S} ]$ may be  singular on $S$. 
But this does not affect the proof, since in the $L^{2}$-extension theorem
(Theorem \ref{l2}) we do not need to assume that the manifold $M$ is compact. 
Hence we may remove a suitable subvarieties so that we do not need 
to consider the pole of $dV[\Psi_{S}]$ on $S$ (but of course the pole of 
$dV[\Psi_{S}]$  affects the $L^{2}$-conditions). 

\begin{theorem}\label{subad2}
Let $M$ be a projective manifold with a continuous volume form $dV$,
let $L$ be a holomorphic line bundle over $M$ with a $C^{\infty}$-hermitian metric  $h_{L}$ with semipositive curvature $\Theta_{h_{L}}$, let $S$ be a compact complex submanifold of $M$,
let $\Psi_{S} : M \longrightarrow [-\infty ,0)$ be a continuous function and let $K_{M}$ be the canonical bundle of $M$.
\begin{enumerate}
\item $\Psi_{S} \in  \sharp (S) \cap C^{\infty}(M\backslash S)$ (As for the 
definition of $\sharp (S)$, see Section 3.2), 
\item $\Theta_{h_{L}\cdot e^{-(1+\epsilon )\Psi_{S}}}\geqq 0$ for 
every $\epsilon \in [0,\delta ]$ for some $\delta > 0$,
\item there is a positive line bundle on $M$.
\end{enumerate}
Then every element of  $H^{0}(S,{\cal O}_{S}(m(K_{M}+L)))$ extends to an element of \\ 
$H^{0}(M,{\cal O}_{M}(m(K_{M}+L)))$. 
$\square$ \end{theorem}
For the completeness we shall give a simple proof of Theorem \ref{subad2} (hence also Theorem \ref{subad1}) 
under the additional conditions :  \vspace{5mm} \\
{\bf Conditions} 
\begin{enumerate}
\item $K_{M}+L$ is big.
\item Bs$\mid m(K_{M}+L)\mid$ does not contain $S$ 
for some $m > 0$.
\item There exists a Zariski open neighbourhood $U$ of the generic
point of $S$ in $M$ such that 
$\mid m(K_{M}+L)\mid$ gives an embedding of $U$ 
into a projective space for every sufficiently large $m$. 
\end{enumerate}
The reason why we put this condition is that we only need Theorems \ref{subad1} and \ref{subad2} under this condition.  More precisely we need to consider the a little bit
more general  case that 
$h_{L}$ is a singular hermitian metric with semipositve curvature current on 
$M$ and $dV[\Psi]$ is singular on $S$. 
But as we have already mentioned above the singularity of $dV[\Psi ]$ 
does not change the proof.  And  the singularity of $h_{L}$ 
will be managed in Remark \ref{r2.6} below. \vspace{5mm} \\

Let us begin the proof of Theorem \ref{subad2} under the above additional conditions. 
Let $M, S, L$ be as in Theorem \ref{subad2}.
Let $n$ denote the dimension of $M$ and let $k$ denote 
the dimension of $S$.   Let $h_{S}$ be a canonical AZD
(\cite{tu2}) of 
$K_{M}+L\mid_{S}$. 
By Kodaira's lemma (cf. \cite[Appendix]{k-o}), there exists an effective {\bf Q}-divisor $B$ on $M$ such that 
$K_{M}+L - B$ is ample.
By the above conditions, we may take $B$ such that 
$\mbox{Supp}\, B$ does not contain $S$. 
In fact  by the conditions, we see that for an ample line bundle 
$H$, $\mid m(K_{M}+L)-H\mid$ is base point free on the generic point 
of $S$. Then we may take $B$ to be the $1/m$-times a general member of 
$\mid m(K_{M}+L)-H\mid$. 
We shall assume that $\mbox{Supp}\, B$ does not contain $S$.

Let $a$ be a positive integer such that 
\begin{enumerate}
\item $A : = a(K_{M}+L-B)$ is Cartier, 
\item  $A\mid_{S}-K_{S}$ ia ample and 
${\cal O}_{S}(A\mid_{S}-K_{S})\otimes {\cal M}_{x}^{k+1}$ is globally generated
for every $x\in S$. 
\end{enumerate}
Let $h_{M}$ be a  canonical AZD of $K_{M}+L$.
We shall define a sequence of the hermitian metrics $\{ \tilde{h}_{m}\} (m\geqq 1)$ inductively by : 
\begin{eqnarray*}
\tilde{K}_{m} &:= & K(M,A+m(K_{M}+L),dV^{-1}\cdot h_{L}\cdot \tilde{h}_{m-1},dV), 
\vspace{5mm} \\
\tilde{h}_{m} & := & \frac{1}{\tilde{K}_{m}},
\end{eqnarray*}
where $K(M,A+m(K_{M}+L),dV^{-1}\cdot h_{L}\cdot \tilde{h}_{m-1},dV)$
is the Bergman kernel of $A+m(K_{M}+L)$ with respect to 
the singular hermitian metric $dV^{-1}\cdot h_{L}\cdot \tilde{h}_{m-1}$ 
and the volume form $dV$, i.e., 
\[
K(M,A+m(K_{M}+L),dV^{-1}\cdot h_{L}\cdot \tilde{h}_{m-1},dV) 
= \sum_{j}\mid\tilde{\sigma}_{j}^{(m)}\mid^{2},
\]
where $\{\tilde{\sigma}^{(m)}_{j}\}$ is a complete orthonormal basis 
of \\ $H^{0}(M,{\cal O}_{M}(A+m(K_{M}+L))\otimes {\cal I}(\tilde{h}_{m-1}))$
with respect to the inner product 
\[
(\tilde{\sigma},\tilde{\sigma}^{\prime}):= \int_{M}\tilde{\sigma}\cdot\bar{\tilde{\sigma}}^{\prime}
\cdot (dV^{-1}\cdot h_{L}\cdot \tilde{h}_{m-1}) \cdot dV
\]
where $\tilde{\sigma},\tilde{\sigma}^{\prime}
\in H^{0}(M,{\cal O}_{M}(A+m(K_{M}+L))\otimes {\cal I}(\tilde{h}_{m-1}))$.
We use the similar notation for Bergman kernels hereafter. 

Every $\tilde{h}_{m}$ is a singular hermitian metric on $A+mK_{M}$ 
with semipositive curvature current by definition.

\begin{lemma}\label{2.2}
For every $m\geqq 0$, there exists a positive constant $C_{m}$ such that 
\[
\tilde{h}_{m}\mid_{S} \leqq C_{m}\cdot h_{A}\mid_{S}\cdot h_{S}^{m}
\]
holds. $\square$  \end{lemma}
{\bf Proof.} 
We shall prove the lemma by induction on $m$. 
For $m= 0$ the both sides are ${\cal O}_{S}$, hence the inclusion holds.
Suppose that the inclusion holds for some $m -1 \geqq 0$ and a positive constant $C_{m-1}$.
Then by the $L^{2}$-extension theorem, Theorem \ref{l2} implies that 
there exists a bounded interpolation operateor : 
\[
I_{m} :  A^{2}(S,A+m(K_{M}+L)\mid_{S},(dV^{-1}\cdot h_{L})\mid_{S}\cdot \tilde{h}_{m-1}\mid_{S},dV[\Psi_{S} ])
\]
\[
\hspace{40mm} 
\longrightarrow 
A^{2}(M,A+ m(K_{M}+L),(dV^{-1}\cdot h_{L})\cdot \tilde{h}_{m-1},dV)
\]
whose operator norm is bounded from above by $C\cdot\delta^{-3/2}$, where $C$ is the 
positive constant in Theorem \ref{l2}. 
Hence by the induction assumption, we see that 
there exists a bounded interpolation operator :
\[
I_{m}^{\prime} :  A^{2}(S,A+m(K_{M}+L)\mid_{S},(dV^{-1}\cdot h_{L})\mid_{S}
\cdot (h_{A}\mid_{S}\cdot h_{S}^{m-1}),dV[\Psi_{S} ])
\]
\[
\hspace{40mm} 
\longrightarrow 
A^{2}(M,A+ m(K_{M}+L),(dV^{-1}\cdot h_{L})\cdot \tilde{h}_{m-1},dV)
\]
whose operator norm is bounded from above by $C_{m-1}\cdot C\cdot\delta^{-3/2}$.
Let $K(S,A+m(K_{M}+L)\mid_{S},(dV^{-1}\cdot h_{L}\cdot h_{A})\mid_{S}\cdot h_{S}^{m-1},dV[\Psi_{S} ])$ denote the Bergman kernel of 
$A+m(K_{M}+L)\mid_{S}$ with respect to the singular hermitian metric 
$(dV^{-1}\cdot h_{L}\cdot h_{A})\mid_{S}\cdot h_{S}^{m-1}$ and 
the volume form $dV[\Psi_{S}]$ (defined as 
$K(M,A+m(K_{M}+L),dV^{-1}\cdot h_{L}\cdot \tilde{h}_{m-1},dV)$ above). 
Then since for every $x\in S$ 
\[
\tilde{K}_{m}(x)
=   \sup\{\mid\tilde{\sigma}(x)\mid^{2} ;
\tilde{\sigma}\in A^{2}(M,A + m(K_{M}+L),dV^{-1}\cdot h_{L}\cdot h_{m-1},dV), \parallel\tilde{\sigma}\parallel = 1\} ,
\]
and 
\[
K(S,A+m(K_{M}+L)\mid_{S},(dV^{-1}\cdot h_{L}\cdot h_{A})\mid_{S}\cdot h_{S}^{m-1},dV[\Psi_{S} ])(x) 
\]
\[
=   \sup\{\mid\sigma (x)\mid^{2} ;  
\sigma\in 
A^{2}(S,A+m(K_{M}+L)\mid_{S},(dV^{-1}\cdot 
h_{L}\cdot h_{A})\mid_{S}\cdot h_{S}^{m-1},dV[\Psi_{S} ]),
 \parallel\sigma\parallel = 1\}
\]
hold (cf. \cite[p.46, Proposition 1.3.16]{kr}), 
we see that there exists a positive constant $C$ such that 
\[
\tilde{K}_{m}\mid_{S} \geqq (C\cdot\delta^{-3/2})^{-1}\cdot C_{m-1}^{-1}\cdot 
K(S,A+m(K_{M}+L)\mid_{S},(dV^{-1}\cdot h_{L}\cdot h_{A})\mid_{S}\cdot h_{S}^{m-1},dV[\Psi_{S} ])
\]
holds on $S$. 
Since there exists a positive constant $C_{1}$ such that 
\[
dV^{-1}\cdot h_{L} \leqq C_{1}\cdot h_{S}
\]
holds, we see that 
\[
\mbox{($\sharp$)}\hspace{5mm}\tilde{K}_{m}\mid_{S} \geqq (C\cdot\delta^{-3/2})^{-1}\cdot C_{m-1}^{-1}\cdot 
C_{1}^{-1}\cdot K(S,A+m(K_{M}+L)\mid_{S},h_{A}\mid_{S}\cdot h_{S}^{m},dV[\Psi_{S} ])
\]
holds. 
By the choice of $A$, we see that there exists a positive constant $C_{S}$ 
(independent of $m$, although this fact is not used in the proof) such that 
\[
\mbox{($\flat$)}\hspace{15mm} K(S,A+m(K_{M}+L)\mid_{S},h_{A}\mid_{S}\cdot h_{S}^{m},dV[\Psi_{S} ])
\geqq C_{S}\cdot (h_{A}\mid_{S}\cdot h_{S}^{m})^{-1}
\]
holds.
This can be verified as follows.  
Since $A\mid_{S}-K_{S}$ is ample, we see that there exists a 
$C^{\infty}$-hermitian metric $h_{A/S}$ on $A\mid_{S}$ such that 
the hermitian metric $dV[\Psi_{S}]\cdot h_{A/S}$
on $A\mid_{S}-K_{S}$ has strictly positive curvature everywhere on $S$.

Let $x$ be a point on $M$ and  $\{\sigma_{A,q}\}$  a basis of 
$H^{0}(S,{\cal O}_{S}(A\mid_{S}-K_{S})\otimes {\cal M}_{x}^{k+1})$.
Then in Theorem \ref{l2} (see also Lemma \ref{l2lemma}), taking $\Psi$ to be 
\[
\Psi_{x} := \frac{k}{k+1}\log \sum_{q}dV[\Psi_{S}]\cdot h_{A/S}(\sigma_{A,q},\sigma_{A,q}),
\]
and $(E,h_{E})$ to be 
\[(A\mid_{S} -K_{S}+ m(K_{M}+L)\mid_{S},dV[\Psi_{S}]\cdot h_{A/S}\cdot h_{S}^{m}),
\]
by Theorem \ref{l2} and Lemma \ref{l2lemma}, we have a bounded interpolation operator :
\[
I_{m,x} : A^{2}(x,A + m(K_{M}+L)\mid_{x},h_{A/S}\cdot h_{S}^{m}\mid_{x},\delta_{x})
\longrightarrow A^{2}(S,A + m(K_{M}+L),h_{A/S}\cdot h_{S}^{m},dV[\Psi_{S}]),
\]
where $\delta_{x}$ is the Dirac measure at $x$.
We note that by the definition of $\Psi_{x}$ and the fact that 
${\cal O}_{S}(A\mid_{S}-K_{S})\otimes {\cal M}_{x}^{k+1}$
is globally generated, $\log\Psi_{x}$ has singularity only at $x$ and  
the operator norm of the $I_{m,x}$ is less than or equal to $C\cdot k^{3/2}$ by Theorem \ref{l2},
where $C$ is the positive constant in Theorem \ref{l2}.
Hence we see that 
\[
K(S,A+m(K_{M}+L)\mid_{S},h_{A/S}\cdot h_{S}^{m},dV[\Psi_{S} ])
\geqq C^{-1}\cdot k^{-3/2}\cdot (h_{A/S}\cdot h_{S}^{m})^{-1}
\]
holds by the basic property of Bergman kernels (cf.\cite[p.46, Proposition 1.3.16]{kr}).
We note that $h_{A/S}$ is quasi-isometric to $h_{A}\mid_{S}$, i.e., 
there exists a positive constant $C_{A,S} > 1$ such that 
\[
C_{A,S}^{-1}\cdot h_{A/S}\leqq   h_{A}\mid_{S} \leqq C_{A,S}\cdot h_{A/S}
\]
holds on $S$. 
Then this implies that 
\[
K(S,A+m(K_{M}+L)\mid_{S},h_{A}\mid_{S}\cdot h_{S}^{m},dV[\Psi_{S} ])
\geqq C_{A,S}^{-1}\cdot C^{-1}\cdot k^{-3/2}\cdot (h_{A}\mid_{S}\cdot h_{S}^{m})^{-1}
\]
holds on $S$. 
This is  the desired estimate ($\flat$) with 
$C_{S} = C_{A,S}^{-1}\cdot C^{-1}\cdot k^{-3/2}$.

Combining ($\sharp$) and ($\flat$), we see that 
\[
\tilde{K}_{m}\mid_{S} \geqq (C\cdot\delta^{-3/2})^{-1}\cdot C_{m-1}^{-1}\cdot 
C_{1}^{-1}\cdot C_{S}\cdot (h_{A}\mid_{S}\cdot h_{S}^{m})^{-1}
\]
holds on $S$. 
Then by the definition of $\tilde{h}_{m}$, we see that 
\[ 
\tilde{h}_{m}\mid_{S} \leqq (C\cdot\delta^{-3/2})\cdot
C_{1}\cdot C_{S}^{-1}\cdot C_{m-1}\cdot  h_{A}\mid_{S}\cdot h_{S}^{m}
\]
holds.
Hence we  complete the proof of Lemma \ref{2.2} by induction on $m$. $\square$ \vspace{5mm}\\

By the definition of $A$, we may consider the metric $h_{A}$ as 
a singular hermitian metric $\hat{h}_{A}$ on $a(K_{M}+L)$. 
Also we may consider $\tilde{h}_{m}$ as a singular hermitian metric 
on $\hat{h}_{m}$ on $(a + m)(K_{M}+L)$.
Then by Lemma \ref{2.2}, we have the following lemma.

\begin{lemma}\label{2.3}
For every $m\geqq 0$, there exist a positive constant $C_{m}^{\prime}$ depending on $m$ and a positive constant $C$ independent of $m$   
such that 
\[
h_{M}^{a+m}\mid_{S}\leqq C_{m}^{\prime}\cdot \hat{h}_{m}\mid_{S} \leqq C^{m+1}\hat{h}_{A}\mid_{S}\cdot h_{S}^{m}
\]
hold. $\square$
\end{lemma}
By Lemma \ref{2.3}, we see that
\[
h_{M}\leqq (C_{m}^{\prime})^{\frac{1}{a+m}}\cdot\hat{h}_{A}\mid_{S}^{\frac{1}{a+m}}\cdot h_{S}^{\frac{m}{a+m}}
\]
holds. 

Let us fix an arbitrary nonnegative integer $\ell$. 
Then  since $h_{S}$ is an AZD of $K_{M}+L\mid_{S}$,
\[
\{{\cal I}(\hat{h}_{A}\mid_{S}^{\frac{\ell}{a+m}}\cdot h_{S}^{\frac{m}{a+m}\ell})\}_{m=1}^{\infty}
\]
is an increasing sequence of ideal sheaves on $S$
contained in ${\cal I}(h_{S}^{\ell})$. 
Let $\phi ,\rho$ be a weight functions of $h_{S}^{\frac{m}{a+m}\ell}$ 
and $\hat{h}_{A}\mid_{S}^{\frac{\ell}{a+m}}$  with respect to (the powers of)
$dV^{-1}\cdot h_{L}\mid_{S}$ respectively. 
By H\"{o}lder's inequality we see that 
for a holomorphic function $f$ on an open set $V$ in $S$, 
\[
\int_{V}e^{-\phi}\cdot e^{-\rho}\cdot\mid f\mid^{2}dV[\Psi_{S} ] \leqq 
(\int_{V} e^{-p\phi}\cdot \mid f\mid^{2}dV[\Psi_{S}])^{\frac{1}{p}}
\cdot (\int_{V} e^{-q\rho}\cdot\mid f\mid^{2}dV[\Psi_{S}])^{\frac{1}{q}}
\]
holds, where 
\[
p := (1+\frac{1}{\ell})(1+\frac{m}{a}), q = \frac{p}{p-1}.
\]
Since 
\[
e^{-p\phi}\cdot (dV^{-1}\cdot h_{L}\mid_{S})^{\ell +1} = h_{S}^{\ell +1}
\]
holds, 
this implies that there exists a positive integer $m_{\ell}$ depending on $\ell$ such that
\[
{\cal I}(\hat{h}_{A}\mid_{S}^{\frac{\ell}{a+m_{\ell}}}\cdot h_{S}^{\frac{m_{\ell}}{a+m_{\ell}}\ell}) 
\supseteq{\cal I}(h_{S}^{\ell +1})
\]
holds. 
Hence we see that 
\[
{\cal I}(h_{M}\mid_{S}^{\ell}) \supseteq  {\cal I}(h_{S}^{\ell +1})
\]
holds on $S$. 
We note that since $h_{S}$ is an AZD of $(K_{M}+L)\mid_{S}$, 
\[
A^{2}(S,(\ell + 1)(K_{M}+L)\mid_{S},h_{S}^{\ell +1},dV[\Psi_{S}])
\simeq 
A^{2}(S,(\ell + 1)(K_{M}+L)\mid_{S},dV^{-1}\cdot h_{L}\mid_{S}\cdot h_{S}^{\ell},dV[\Psi_{S}])
\]
holds. 
Using this equality, by Theorem \ref{l2} (and Lemma \ref{l2lemma}) in Section 2.3, we see that every element of 
\[
A^{2}(S,(\ell + 1)(K_{M}+L)\mid_{S},dV^{-1}\cdot h_{L}\mid_{S}\cdot h_{S}^{\ell},dV[\Psi_{S}])
\]
can be extended to an element of 
\[
A^{2}(M,(\ell + 1)(K_{M}+L),dV^{-1}\cdot h_{L}\cdot h_{M}^{\ell},dV).
\]
Since $\ell$ is an arbitrary nonnegative integer, we complete 
the proof of Theorem \ref{subad2}. $\square$. 

\begin{remark}\label{r2.6}
The above proof also works for the case that $(L,h_{L})$ is a singular hermitian line bundle with semipositive curvature current, if we assume the following 
conditions : 
\begin{enumerate}
\item $(K_{M}+L,dV^{-1}\cdot h_{L})$ is big.
\item Bs$\mid m(K_{M}+L,dV^{-1}\cdot h_{L})\mid_{\infty}$ does not contain $S$ 
for some $m > 0$.
\item There exists a Zariski open neighbourhood $U$ of the generic
point of $S$ in $M$ such that 
$\mid m(K_{M}+L,dV^{-1}\cdot h_{L})\mid_{\infty}$ gives an embedding of $U$ 
into a projective space for every sufficiently large $m$.
\end{enumerate}
Here $\mid m(K_{M}+L,dV^{-1}\cdot h_{L})\mid_{\infty}$ denotes the linear system 
$\mid H^{0}(M,{\cal O}_{M}(m(K_{M}+L))\otimes {\cal I}_{\infty}(h_{L}^{m}))\mid$. 
In this case we need to take an AZD $h_{S}$ of the singular hermitian line bundle 
$(K_{M}+L,dV^{-1}\cdot h_{L})\mid_{S}$. 
Noting Remarks \ref{r2.3} and \ref{r2.4}, by Corollary \ref{cor} there exists an 
AZD  $h_{S}$ of $(K_{M}+L,dV^{-1}\cdot h_{L})\mid_{S}$. $\square$
\end{remark}
\begin{remark}
The full proofs of Theorems \ref{subad1} and \ref{subad2} can be obtained 
similar line as the above proof. 
But they  require more detailed estimates. 
The proof presented here is somewhat similar to the argument in \cite{si}. $\square$ 
\end{remark}
\subsection{Positivity result}

The following positivity theorem is the key to the proof of Theorems 1.1 and 
1.2.
\begin{theorem}\label{pos}(\cite[p.894,Theorem 2]{ka})
Let $f : X \longrightarrow B$ be a surjective morphism of smooth projective 
varieties with connected fibers.
Let $P = \sum P_{j}$ and $Q = \sum_{\ell}Q_{\ell}$ be normal crossing divisors on $X$ and $B$ respectively, such that $f^{-1}(Q) \subset P$ and $f$ 
is smooth over $B\backslash Q$.
Let $D = \sum d_{j}P_{j}$ be a {\bf Q}-divisor on $X$, where $d_{j}$ may be positive, zero or negative, which satisfies the following conditions :
\begin{enumerate}
\item $D = D^{h} + D^{v}$ such that 
$f :\mbox{Supp}(D^{h})\rightarrow B$ is surjective and smooth over $B\backslash Q$, and $f(\mbox{Supp}(D^{v}))\subset Q$.
An irreducible component of $D^{h}$(resp. $D^{v}$) is called horizontal
(resp. vertical).
\item $d_{j} < 1$ for all $j$.
\item The natural homomorphism ${\cal O}_{B}\rightarrow f_{*}{\cal O}_{X}(\lceil -D\rceil )$ is surjective at the generic point of $B$.
\item $K_{X} +  D\sim_{\mbox{\bf Q}}f^{*}(K_{B} + L)$ for some 
{\bf Q}-divisor $L$ on $B$.
\end{enumerate} 
Let 
\begin{eqnarray*}
f^{*}Q_{\ell}& =  &\sum_{j}w_{\ell j}P_{j} \\
\bar{d}_{j} & :=  & \frac{d_{j} +w_{\ell j}-1}{w_{\ell j}}\,\,\,\,\mbox{if}\,\,\,\,
f(P_{j}) = Q_{\ell} \\
\delta_{\ell} &: =  & \max \{\bar{d}_{j} ; f(P_{j}) = Q_{\ell}\} \\
\Delta & :=  & \sum_{\ell}\delta_{\ell}Q_{\ell} \\
M & :=  & L - \Delta .
\end{eqnarray*}
Then $M$ is nef.  $\square$

 \end{theorem} 
 
Here the meaning of the divisor $\Delta$ may be difficult to understand.
So I would like to give an geometric interpretation of $\Delta$.  
Let $X,P,Q,D,B,\Delta$ be as above. Let $dV$ be a 
$C^{\infty}$-volume form on $X$. 
Let $\sigma_{j}$ be a global section of ${\cal O}_{X}(P_{j})$
with divisor $P_{j}$. 
Let $\parallel\sigma_{j}\parallel$ denote the hermitian norm 
of $\sigma_{j}$ with respect to a $C^{\infty}$-hermitian metric
on ${\cal O}_{X}(P_{j})$ respectively. 
Let us consider the singular volume form
\[
\Omega := \frac{dV}{\prod_{j}\parallel\sigma_{j}\parallel^{2d_{j}}}
\]
on $X$.
Then by taking the fiber integral of $\Omega$ with respect to 
$f : X \longrightarrow B$, we obtain a singular volume form 
$\int_{X/B}\Omega$ on $B$, where the fiber integral $\int_{X/B}\Omega$
is defined by the property that for any open set $U$ in $B$, 
\[
\int_{U}(\int_{X/B}\Omega ) = \int_{f^{-1}(U)}\Omega
\]
holds. 
We note that the  condition 2 in Theorem \ref{pos} assures that 
$\int_{X/B}\Omega$ is continuous on a nonempty Zariski open subset 
of $B$.
Also by the condition 4 in Theorem \ref{pos}, we see that 
$K_{X}+D$ is numerically $f$-trivial and 
$(\int_{X/B}\Omega)^{-a}$ is  a $C^{0}$-hermitian metric 
on a line bundle $a(K_{B}+\Delta )$, where $a$ is a positive integer 
such that $a\Delta$ is Cartier. 
Thus  the divisor $\Delta$ corresponds exactly to  
singularities (poles and degenerations)  
of the singular volume form $\int_{X/B}\Omega$ on $B$.

\section{Proofs of Theorems 1.1 and 1.2}

In this section we shall prove Theorems 1.1 and 1.2
simultaneously. The proof is  almost parallel to the one  assuming 
MMP (\cite{tu8}), if we replace the minimal model by an AZD (analytic
Zariski decomposition) of the canonical line bundle. 

\subsection{Construction of a filtration}

Let $X$ be a smooth projective $n$-fold of general type.
Let $h$ be an AZD of $K_{X}$ constructed as in Section 2.2. 
We may assume that $h$ is lowersemicontinuous by Theorem \ref{lelong}. 
This is a technical assumption so that a local potential 
of the curvature current of $h$ is plurisubharmonic. 
This is used to restrict $h$ to a subvariety of $X$ 
(if we only assume that the  local potential is only locally integrable, 
the restriction is not well defined). 
We set 
\[
X^{\circ} = \{ x\in X\mid x\not{\in} \mbox{Bs}\mid mK_{X}\mid \mbox{and  
$\Phi_{\mid mK_{X}\mid}$ is a biholomorphism} 
\]
\[
\hspace{50mm} \mbox{on a neighbourhood of $x$ for some $m \geqq 1$}\} .
\]
We set 
\[
\mu_{0} := (K_{X},h)^{n}= \mu (X,(K_{X},h)) = \mu (X,K_{X}).
\]
For the notations $(K_{X},h)^{n}, \mu (X,(K_{X},h))$ and $\mu (X,K_{X})$ 
see Definitions \ref{volume} and \ref{volume1}. 
The last equality holds, since $h$ is an AZD of $K_{X}$.
We note that for every $x\in X^{\circ}$, 
${\cal I}(h^{m})_{x}\simeq {\cal O}_{X,x}$ holds for every $m\geqq 0$ 
(cf. \cite{tu2} or \cite[Theorem 1.5]{d-p-s}).  \vspace{5mm}\\
Let $x,x^{\prime}$ be distinct points on $X^{\circ}$. 
In this subsection we shall construct a filtration 
\[
X = X_{0}\supset X_{1}\supset \cdots \supset X_{r}\supset X_{r+1} = 
 x\,\,\mbox{or}\,\, x^{\prime}
\]
of $X$ 
by a strictly decreasing sequence of subvarieties $\{ X_{i}\}_{i=0}^{r+1}$
for some $r$ (depending $x$ and $x^{\prime}$)
and invariants :
\[
\alpha_{0} ,\alpha_{1},\ldots ,\alpha_{r},
\]
\[
n =: n_{0} >  n_{1}> \cdots > n_{r} \,\,\,\,\,(n_{i} = \dim X_{i},i=0,\cdots ,r)
\]
and
\[
\mu_{0},\mu_{1},\ldots ,\mu_{r} \,\,\,\,\,(\mu_{i}= (K_{X},h)^{n_{i}}\cdot X_{i}, i= 0,\cdots ,r)
\]
(cf. Definition \ref{volume}) with the estimates
\[
\alpha_{i} \leqq \frac{n_{i}\sqrt[n_{i}]{2}}{\sqrt[n_{i}]{\mu_{i}}} + \delta
\,\,\,\,\,\, (0\leqq i\leqq r),
\]
where $\delta$ is a fixed positive number less than $1/n$. 

\begin{lemma}  
We set 
\[
{\cal M}_{x,x^{\prime}} = {\cal M}_{x}\cdot{\cal M}_{x^{\prime}},
\]
 where ${\cal M}_{x},{\cal M}_{x^{\prime}}$ denote the
maximal ideal sheaves of the points $x,x^{\prime}$ respectively.
Let $\varepsilon$ be  positive number strictly less than $1$.
Then 
\[
H^{0}(X,{\cal O}_{X}(mK_{X})\otimes {\cal I}(h^{m})\cdot {\cal M}_{x,x^{\prime}}^{\lceil\sqrt[n]{\mu_{0}}
(1-\varepsilon )\frac{m}{\sqrt[n]{2}}\rceil})\neq 0
\]
holds for every sufficiently large $m$.
$\square$ \end{lemma}
{\bf Proof}.
First we note that since $x,x^{\prime}\in X^{\circ}$, 
$h$ is bounded from above at $x$ and $x^{\prime}$ by the construction of $h$
(cf. Theorem \ref{AZD}).  In particular ${\cal I}(h^{m})_{x} = {\cal O}_{X,x}$ 
and ${\cal I}(h^{m})_{x^{\prime}} = {\cal O}_{X,x^{\prime}}$ hold 
for every $m \geqq 0$. 
Let us consider the exact sequence:
\[
0\rightarrow H^{0}(X,{\cal O}_{X}(mK_{X})\otimes {\cal I}(h^{m})\cdot 
{\cal M}_{x,x^{\prime}}^{\lceil\sqrt[n]{\mu_{0}}(1-\varepsilon )\frac{m}{\sqrt[n]{2}}\rceil})
\rightarrow H^{0}(X,{\cal O}_{X}(mK_{X})\otimes {\cal I}(h^{m}))\rightarrow
\]
\[
  H^{0}(X,{\cal O}_{X}
(mK_{X})\otimes {\cal I}(h^{m})\otimes {\cal O}_{X}/{\cal M}_{x,x^{\prime}}^{\lceil\sqrt[n]{\mu_{0}}(1-\varepsilon )\frac{m}{\sqrt[n]{2}}\rceil}).
\]
We note that 
\[
n!\cdot\overline{\lim}_{m\rightarrow\infty}m^{-n}\dim H^{0}(X,{\cal O}_{X}(mK_{X})\otimes {\cal I}(h^{m})) = \mu_{0}
\]
holds by the definition of $\mu_{0}$.

On the other hand, we see that 
\[
n!\cdot\overline{\lim}_{m\rightarrow\infty}m^{-n}\dim H^{0}(X,{\cal O}_{X}(mK_{X})
\otimes {\cal I}(h^{m})\otimes {\cal O}_{X}/{\cal M}_{x,x^{\prime}}^{\lceil\sqrt[n]{\mu_{0}}(1-\varepsilon )\frac{m}{\sqrt[n]{2}}\rceil})
=
\mu_{0}(1-\varepsilon )^{n} < \mu_{0}
\]
hold, since ${\cal I}(h^{m})_{x} = {\cal O}_{X,x}$ 
and ${\cal I}(h^{m})_{x^{\prime}} = {\cal O}_{X,x^{\prime}}$ hold 
for every $m \geqq 0$.  

By the above inequalities and  the  exact sequence, we complete the proof of  Lemma 3.1.  $\square$ \vspace{5mm}\\
Let $\varepsilon$ be a positive number less than $1$ as in Lemma 3.1.
Let us take a sufficiently large positive integer $m_{0}$ such that 
\[
H^{0}(X,{\cal O}_{X}(m_{0}K_{X})\otimes {\cal I}(h^{m_{0}})\cdot 
{\cal M}_{x,x^{\prime}}^{\lceil\sqrt[n]{\mu_{0}}(1-\varepsilon )\frac{m_{0}}{\sqrt[n]{2}}\rceil}) \neq 0
\]
as in Lemma 3.1 and let 
\[
\sigma_{0} \in H^{0}(X,{\cal O}_{X}(m_{0}K_{X})\otimes {\cal I}(h^{m_{0}})
\cdot{\cal M}_{x,x^{\prime}}^{\lceil\sqrt[n]{\mu_{0}}(1-\varepsilon )\frac{m_{0}}{\sqrt[n]{2}}\rceil})
\]
be a general nonzero element.
We set 
\[
D_{0} := \frac{1}{m_{0}}(\sigma_{0})
\]
and  
\[
h_{0} := \frac{1}{\mid \sigma_{0}\mid^{2/m_{0}}}
\]
(see  Example \ref{ex} in Section 2.1 for the meaning of $1/\mid \sigma_{0}\mid^{2/m_{0}}$ ).
We define the positive number $\alpha_{0}$ by
\[
\alpha_{0} := \inf\{\alpha > 0\mid 
\mbox{$(X,\alpha D_{0})$ is KLT at 
neither $x$ nor $x^{\prime}$}\},
\]
where KLT is short for of Kawamata log terminal (cf. Definition \ref{KLT}).  

Let $\mu : Y \longrightarrow X$ be a log resolution of $(X,D)$ 
and for $\alpha > 0$ let 
\[
K_{Y} + \mu_{*}^{-1}(\alpha D) = \mu^{*}(K_{X}+\alpha D) + F(\alpha ),
\]
where $F(\alpha )$ denotes the discrepancy depending on $\alpha$. 
Then $\alpha_{0}$ is the infimum of $\alpha$ such that 
the discrepancy $F(\alpha )$ has a component whose coefficient is less than or 
equal to $-1$. 
Hence by the construction $\alpha_{0}$ is a rational number. 

Since $(\sum_{i=1}^{n}\mid z_{i}\mid^{2})^{-n}$ is not locally integrable 
around $O\in \mbox{{\bf C}}^{n}$, by the definition of $D_{0}$, we see
that 
\[
\alpha_{0}\leqq \frac{n\sqrt[n]{2}}{\sqrt[n]{\mu_{0}}(1-\varepsilon )}
\]
holds.  About the relation between the KLT conditions and 
the multiplier ideal sheaves, please see Section 2.1.

Let $\delta$ be the fixed positive number as above 
and let us make  $\varepsilon > 0$ sufficiently small so that 
\[
\alpha_{0} \leqq \frac{n\sqrt[n]{2}}{\sqrt[n]{\mu_{0}}}+ \delta 
\] 
holds. 
Then one of the following two cases occurs. \vspace{5mm} \\
{\bf Case} 1:  For every sufficiently small positive number $\lambda$, 
$(X,(\alpha_{0}-\lambda )D_{0})$ is KLT 
at both $x$ and $x^{\prime}$. \\
{\bf Case} 2: For every sufficiently small positive number $\lambda$, 
$(X,(\alpha_{0}-\lambda )D_{0})$ is KLT  
at exactly one of $x$ or $x^{\prime}$, say $x$. \vspace{5mm} \\
We define the next stratum $X_{1}$ by 
\[
X_{1}:=  \mbox{the minimal center of log canonical singularities 
of $(X,\alpha_{0}D_{0})$ at $x$  (cf. Section 2.1).}
\] 
If $X_{1}$ is a point, we stop the construction of the filtration. 
Suppose that $X_{1}$ is not a point. 

Case 1 divides into  the following two cases. \vspace{5mm} \\
{\bf Case} 1.1: $X_{1}$ passes through both $x$ and $x^{\prime}$, \\
{\bf Case} 1.2: $X_{1}$ passes through only one of $x$ and $x^{\prime}$, say $x$. \vspace{5mm} \\
First we shall consider Case1.1.
We define the positive number $\mu_{1}$ by 
\[
\mu_{1}:= \mu (X_{1},(K_{X},h)\mid_{X_{1}}).
\]
Then since $x,x^{\prime}\in X^{\circ}$, 
$\mu_{1}$ is positive. 

For the later purpose, we shall modify  $h_{0}$ 
so that $X_{1}$ is the only center of log canonical singularities 
of $(X,\alpha_{0}D_{0})$ at $x$. 
Let us take an effective {\bf Q}-divisor $G$ such that 
$K_{X} - G$ is ample by Kodaira's lemma (cf.\cite[Appendix]{k-o}). 
By the definition of $X^{\circ}$, we may assume that the support of $G$ contains neither $x$ nor $x^{\prime}$. 
In fact this can be verified as follows. Let $H$ be an arbitrary ample 
divisor on $X$. 
Then by the definition of $X^{\circ}$, $\mid bK_{X} - H\mid$ is base point free at $x$ and $x^{\prime}$ 
for every sufficiently large $b$.
Fix such a $b$ and take a member $G^{\prime}$ of $\mid b K_{X} - H\mid$
which contains neither $x$ nor $x^{\prime}$. 
Then we may take $G$ to be  $b^{-1}G^{\prime}$.

Let $a$ be a positive integer such that 
$A:= a (K_{X}- G)$ is a very ample Cartier divisor 
such that ${\cal O}_{X}(A)\otimes {\cal I}_{X_{1}}$ is globally 
generated. 
Let $\rho_{1}, \ldots ,\rho_{e} 
\in H^{0}(X,{\cal O}_{X}(A)\otimes {\cal I}_{X_{1}})$ be a set of generators 
of ${\cal O}_{X}(A)\otimes {\cal I}_{X_{1}}$ on $X$. 
Then if we replace $h_{0}$ by 
\[
\frac{1}{(\mid \sigma_{0}\mid^{2}(\sum_{i=1}^{e}\mid\rho_{i}\mid^{2}))^{\frac{1}{m_{0}+a}}}, 
\]
it has the desired property. 
If we take $m_{0}$ very large (in comparison with $a$), we can make 
the new $\alpha_{0}$ arbitrary close to the original $\alpha_{0}$.
Hence we may assume that the estimate 
\[
\alpha_{0}\leqq \frac{n\sqrt[n]{2}}{\sqrt[n]{\mu_{0}}}+\delta
\]
still holds after the modification. 
Let us set
\[
n_{1}:= \dim X_{1}.
\]
The proof of the following lemma is identical to that of Lemma 3.1.
\begin{lemma}
Let $\varepsilon^{\prime}$ be a positive rational number less than $1$.
Let $x_{1},x_{2}$ be distinct regular points  of $X_{1}\cap X^{\circ}$. 
Then for every  sufficiently large positive integer  $m$
\[
H^{0}(X_{1},{\cal O}_{X_{1}}(mK_{X})\otimes {\cal I}(h^{m}\mid_{X_{1}})
\cdot{\cal M}_{x_{1},x_{2}}^{\lceil\sqrt[n_{1}]{\mu_{1}}(1-\varepsilon^{\prime} )\frac{m}{\sqrt[n_{1}]{2}}\rceil}) \neq 0
\]
holds.
$\square$ \end{lemma}
Let $x_{1},x_{2}$ be two distinct regular points of  $X_{1}\cap X^{\circ}$. 
Let $m_{1}$ be a positive integer 
such that
\[
H^{0}(X_{1},{\cal O}_{X_{1}}(m_{1}K_{X})\otimes {\cal I}(h^{m_{1}}\mid_{X_{1}})
\cdot{\cal M}_{x_{1},x_{2}}^{\lceil\sqrt[n_{1}]{\mu_{1}}(1-\varepsilon^{\prime} )\frac{m}{\sqrt[n_{1}]{2}}\rceil}) \neq 0
\]
as in Lemma 3.2 
and let 
\[
\sigma_{1,x_{1},x_{2}}^{\prime}
\in 
H^{0}(X_{1},{\cal O}_{X_{1}}(m_{1}K_{X})\otimes {\cal I}(h^{m_{1}}\mid_{X_{1}})
\cdot{\cal M}_{x_{1},x_{2}}^{\lceil\sqrt[n_{1}]{\mu_{1}}(1-\varepsilon^{\prime} )\frac{m}{\sqrt[n_{1}]{2}}\rceil})
\]
be a nonzero element.  

We shall extend the singular hermitian metric 
$1/\mid \sigma_{1,x_{1},x_{2}}^{\prime}\mid^{2/m_{1}}$ of $K_{X}\mid_{X_{1}}$ to 
a singular hermitian metric on $K_{X}$ with semipositive curvature current
after a modification. 

As before by Kodaira's lemma (\cite[Appendix]{k-o}) there is an effective {\bf Q}-divisor $G$ such
that $K_{X}- G$ is ample.
By the definition of $X^{\circ}$, we may assume that 
the support of $G$ contains neither $x$ nor $x^{\prime}$ as before. 
Let $\ell_{1}$ be a sufficiently large positive integer which will be specified later  such that
\[
L_{1} := \ell_{1}(K_{X}- G)
\]
is Cartier.  Let $h_{L_{1}}$ be a $C^{\infty}$-hermitian metric on $L_{1}$ 
with strictly positive curvature. 
Let $\tau$ be a nonzero section in 
$H^{0}(X,{\cal O}_{X}(L_{1}))$.  
We set 
\[
\Psi := \alpha_{0}\cdot\log\frac{h}{h_{0}}.
\]
Let $dV$ be a $C^{\infty}$-volume form on $X$.
We note that the residue volume form  $dV [\Psi ]$ on $X_{1}$
 may have poles along some proper subvarieties in $X_{1}$. 
By taking $\ell_{1}$ sufficiently large and taking $\tau$ properly,
we may assume that $h_{L_{1}}(\tau ,\tau )\cdot dV[\Psi ]$ is nonsingular on 
$X_{1}$ in the sense that the pullback of  it to a nonsingular model of $X_{1}$
is a bounded form. 
Then by applying  Lemma \ref{l2lemma} for $S = X_{1}$, $\Psi = \alpha_{0}\log (h/h_{0})$,
\[
(E,h_{E})= ((\lceil 1+\alpha_{0}\rceil )K_{X},h^{\lceil 1+\alpha_{0}\rceil}),
\]
and
\[
(L,h_{L}) = ((m_{1}-\lceil \alpha_{0}\rceil -2)K_{X}+L_{1},h^{(m_{1}-\lceil \alpha_{0}\rceil -2)}\otimes h_{L_{1}}),
\]
we see that   
\[
\sigma_{1,x_{1},x_{2}}^{\prime}\otimes\tau\in
H^{0}(X_{1},{\cal O}_{X_{1}}(m_{1}K_{X}+L_{1})
\otimes {\cal I}(h^{m_{1}}\mid_{X_{1}})\cdot{\cal M}_{x_{1},x_{2}}^{\lceil\sqrt[n_{1}]{\mu_{1}}(1-\varepsilon^{\prime})\frac{m_{1}}{\sqrt[n_{1}]{2}}\rceil})
\]
extends to a section
\[
\sigma_{1,x_{1},x_{2}}\in H^{0}(X,{\cal O}_{X}((m_{1}+\ell_{1} )K_{X})).
\]
We note that even though $dV[\Psi ]$ may be singular on $X_{1}$, 
we may apply
Lemma \ref{l2lemma}, because there exists a proper Zariski closed subset  
$B$ of $X$ such that the restriction of $dV[\Psi ]$
to $(X - B)\cap X_{1}$ is smooth. 
Of course the singularity of $dV[\Psi ]$ affects  the $L^{2}$-condition.
But this has already been managed by the boundedness of 
$h_{L_{1}}(\tau ,\tau )\cdot dV[\Psi ]$. 

Taking  $\ell_{1}$
sufficiently large, we may and do assume that  there exists a neighbourhood $U_{x_{1},x_{2}}$ of 
$\{ x_{1},x_{2}\}$ such that the divisor $(\sigma _{1,x_{1}.x_{2}})$  is smooth
on  $U_{x_{1},x_{2}} \backslash X_{1}$. 
This can be verified as follows. 
Let us take $\ell_{1}$ sufficiently large so that  ${\cal O}_{X}(L_{1})
\otimes {\cal M}_{y}^{n+1}$ is globally generated for every $y\in X$. 
Let us fix $y$ and let $\{\xi_{1},\cdots ,\xi_{N}\}$ be 
a set of basis of $H^{0}(X,{\cal O}_{X}(L_{1})\otimes {\cal M}_{y}^{n+1})$.
Then 
\[
h_{L_{1},y}:= h_{L_{1}}^{\frac{1}{n+1}}\cdot(\frac{1}{\sum_{j=1}^{N}\mid\xi_{j}\mid^{2}})^{\frac{n}{n+1}}
\]
is a singular hermitian metric of $L_{1}$ with strictly positive curvature 
current. 
Since ${\cal O}_{X}(L_{1})
\otimes {\cal M}_{y}^{n+1}$ is globally generated, we see that 
${\cal O}_{X}/{\cal I}(h_{L_{1},y})$ has isolated support at $y$.
By Nadel's vanishing theorem \cite[p.561]{n}, this implies that 
for every $y\in X^{\circ}\backslash X_{1}$, 
\[
H^{1}(X,{\cal O}_{X}(mK_{X}+L_{1})\otimes {\cal I}(h_{0}^{\alpha_{0}}\cdot
h^{m-1-\alpha_{0}}) \otimes {\cal M}_{y}) = 0
\]
holds. 
Hence   for every $y\in X^{\circ} \backslash X_{1}$, we may modify 
the $L^{2}$-extension of $\sigma_{1,x_{1},x_{2}}^{\prime}\otimes \tau$  so that the extension has any prescribed 
value at $y$, if we take $\ell_{1}$ is sufficiently large.
  We may take $\ell_{1}$ to be independent of 
$y\in X^{\circ} \backslash X_{1}$. 
Then by Bertini's theorem we may find a neighbourhood $U_{x_{1},x_{2}}$ of 
$\{ x_{1},x_{2}\}$ such that the divisor $(\sigma _{1,x_{1}.x_{2}})$  is smooth
on  $U_{x_{1},x_{2}} \backslash X_{1}$.

We set 
\[
D_{1}(x_{1},x_{2}):=  \frac{1}{m_{1}+\ell_{1}}(\sigma_{1,x_{1},x_{2}}). 
\]
Let $X_{1,\mbox{reg}}$ denote the set of regular points on $X_{1}$. 
We may construct the divisors $\{ D_{1}(x_{1},x_{2})\}$ as an algebraic 
family over $(X_{1,\mbox{reg}}\times X_{1,\mbox{reg}})\,\backslash\,\,\Delta_{X_{1}}$ where 
$\Delta_{X_{1}}$ denotes the diagonal of $X_{1}\times X_{1}$.
This construction is possible, since  we may take $L_{1}$  idependent of $x_{1},x_{2}\in X_{1,\mbox{reg}}$. 
Letting $x_{1}$ and $x_{2}$ tend to $x$ and $x^{\prime}$ respectively,
we obtain a {\bf Q}-divisor $D_{1}$ on $X$ which is $(m_{1}+\ell_{1})^{-1}$-times a divisor of a global holomorphic 
section
\[
\sigma_{1}\in H^{0}(X,{\cal O}_{X}((m_{1}+\ell_{1})K_{X})).
\]
By the construction, we may and do assume that there exists a neighbourhood $U_{x,x^{\prime}}$ 
 of $\{ x,x^{\prime}\}$ such that  $(\sigma_{1})$ is smooth 
 on $U_{x,x^{\prime}}\backslash X_{1}$.

Let $\varepsilon_{0}$ be a
positive rational number with $\varepsilon_{0} < \alpha_{0}$.
And we define the positive numbers $\alpha_{1}(x_{1},x_{2})$ and $\alpha_{1}$ by
\[
\alpha_{1}(x_{1},x_{2}) := \inf\{\alpha > 0 \mid 
\mbox{ $(\alpha_{0}-\varepsilon_{0})D_{0} + \alpha D_{1}(x_{1},x_{2})$ 
is  KLT at neither $x_{1}$ nor $x_{2}$} \} 
\]
and
\[
\alpha_{1} := \inf\{\alpha > 0 \mid 
\mbox{ $(\alpha_{0}-\varepsilon_{0})D_{0} + \alpha D_{1}$ 
is  KLT at neither $x$ nor $x^{\prime}$} \}
\]
respectively.
For every positive number 
$\lambda$, $(\alpha_{0}-\varepsilon_{0})D_{0} + (\alpha_{1}-\lambda )D_{1}$ is 
KLT at $x$ or $x^{\prime}$, say $x$. 
Then we shall define the proper subvariety $X_{2}$ of $X_{1}$ by
\[
X_{2}:= \mbox{the minimal center of log canonical singularities 
of $(X,(\alpha_{0}-\varepsilon_{0})D_{0} + \alpha_{1} D_{1})$
at $x$}.
\]

We shall estimate $\alpha_{1}$.
We note that $m_{1}$ is  independent of $\ell_{1}$ 
in the extension of $\sigma_{1,x_{1},x_{2}}^{\prime}\otimes \tau$.

\begin{lemma}
Let $\delta$ be the fixed positive number as above, then we may assume that 
\[
\alpha_{1}\leqq \frac{n_{1}\sqrt[n_{1}]{2}}{\sqrt[n_{1}]{\mu_{1}}} 
+ \delta 
\]
holds, if we take $\varepsilon^{\prime}$, $\ell_{1}/m_{1}$ and $\varepsilon_{0}$ sufficiently small. 
$\square$ \end{lemma}
To prove Lemma 3.3, we need the following elementary lemma.
\begin{lemma}(\cite[p.12, Lemma 6]{t})
Let $a,b$ be  positive numbers. Then
\[
\int_{0}^{1}\frac{r_{2}^{2n_{1}-1}}{(r_{1}^{2}+r_{2}^{2a})^{b}}
dr_{2}
=
r_{1}^{\frac{2n_{1}}{a}-2b}\int_{0}^{r_{1}^{-{2}{a}}}
\frac{r_{3}^{2n_{1}-1}}{(1 + r_{3}^{2a})^{b}}dr_{3}
\]
holds, where 
\[
r_{3} = r_{2}/r_{1}^{1/a}.
\]
$\square$ \end{lemma}
{\bf Proof of Lemma 3.3.}
First suppose that $x,x^{\prime}$ are {\bf nonsingular points} on $X_{1}$.
Then we may set $x_{1} = x, x_{2} = x^{\prime}$, i.e., we do not need the 
limiting process to define the divisor $D_{1}$. 
Let $(z_{1},\ldots ,z_{n})$ be a local coordinate  system on a 
neighbourhood $U$ of $x$ in $X$ such that 
\[
U \cap X_{1} = 
\{ q\in U\mid z_{n_{1}+1}(q) =\cdots = z_{n}(q)=0\} .
\] 
We set $r_{1} = (\sum_{i=n_{1}+1}^{n}\mid z_{1}\mid^{2})^{1/2}$ and 
$r_{2} = (\sum_{i=1}^{n_{1}}\mid z_{i}\mid^{2})^{1/2}$.
Fix an arbitrary $C^{\infty}$-hermitian metric $h_{K}$ on $K_{X}$. 
Then there exists a positive constant $C$ such that 
\[
\mbox{($\star$)}\hspace{10mm} \parallel\sigma_{1}\parallel^{2}\leqq 
C(r_{1}^{2}+r_{2}^{2\lceil\sqrt[n_{1}]{\mu_{1}}(1-\varepsilon^{\prime})\frac{m_{1}}
{\sqrt[n_{1}]{2}}\rceil})
\]
holds on a neighbourhood of $x$, 
where $\parallel\,\,\,\,\parallel$ denotes the norm with 
respect to $h_{K}^{m_{1}+\ell_{1}}$.
We note that there exists a positive integer $M$ such that 
\[
\parallel\sigma_{1}\parallel^{-2} = O(r_{1}^{-M})
\]
holds on a neighbourhood of the generic point of $U\cap X_{1}$,
where $\parallel\,\,\,\,\parallel$ denotes the norm with respect to 
$h_{K}^{m_{0}}$. 
Let us apply Lemma 3.4 by taking 
\[
a = \lceil\sqrt[n_{1}]{\mu_{1}}(1-\varepsilon^{\prime})\frac{m_{1}}
{\sqrt[n_{1}]{2}}\rceil .
\]
Then by Lemma 3.4 and the estimate ($\star$), we see that for every
\[
b > \frac{n_{1}}
{\lceil\sqrt[n_{1}]{\mu_{1}}
(1-\varepsilon^{\prime})\frac{m_{1}}{\sqrt[n_{1}]{2}}\rceil }.
\]
$\parallel\sigma_{1}\parallel$ produces a singularity greater than equal
to $r_{1}^{\frac{2n_{1}}{a} -b}$, if we average the singularity 
in terms of the volume form in $z_{1},\cdots ,z_{n_{1}}$ direction.
Hence by Proposition \ref{prop},  we have the inequality:
\[
\alpha_{1} \leqq (\frac{m_{1}+\ell_{1}}{m_{1}})\frac{n_{1}\sqrt[n_{1}]{2}}{\sqrt[n_{1}]{\mu_{1}}(1-\varepsilon^{\prime})} + m_{1}\varepsilon_{0}.
\]  
Taking $\varepsilon^{\prime}$, $\ell_{1}/m_{1}$ and $\varepsilon_{0}$
sufficiently small, 
we obtain that 
\[
\alpha_{1}\leqq \frac{n_{1}\sqrt[n_{1}]{2}}{\sqrt[n_{1}]{\mu_{1}}} 
+ \delta 
\]
holds.
$\square$ \vspace{5mm} \\
If $x$ or $x^{\prime}$ is a singular point on $X_{1}$, we need the following lemma.
\begin{lemma}
Let $\varphi$ be a plurisubharmonic function on $\Delta^{n}\times{\Delta}$.
Let $\varphi_{t}(t\in\Delta )$ be the restriction of $\varphi$ on
$\Delta^{n}\times\{ t\}$.
Assume that $e^{-\varphi_{t}}$ does not belong to $L^{1}_{\mbox{loc}}(\Delta^{n},O)$
for any $t\in \Delta^{*}$.

Then $e^{-\varphi_{0}}$ is not locally integrable at $O\in\Delta^{n}$.
$\square$ \end{lemma}
Lemma 3.5 is an immediate consequence of
the $L^{2}$-extension theorem \cite[p.20, Theorem]{o-t}.
Using Lemma 3.5 and Lemma 3.4, letting $x_{1}\rightarrow x$ and $x_{2}\rightarrow x^{\prime}$, we see that 
\[
\alpha_{1}\leqq \liminf_{x_{1}\rightarrow x,x_{2}\rightarrow x^{\prime}} \alpha_{1}(x_{1},x_{2})
\]
holds. 
\vspace{5mm}

Next we consider Case 2.
The remaining case Case 1.2 will be considered later.   
In Case 2,   for every  sufficiently small positive number $\lambda$, 
$(X,(\alpha_{0}-\lambda )D_{0})$
is KLT at $x$ and not KLT at $x^{\prime}$.
In Case 1.2, instead of Lemma 3.2, we use the following simpler lemma.
We define $X_{1}$ as before. 
  
In this case, instead of Lemma 3.2, we use the following simpler lemma.

\begin{lemma} Let $\varepsilon^{\prime}$ be a positive number less than $1$ and let $x_{1}$ be a smooth point on $X_{1}$. 
Then for a sufficiently large $m >1$,
\[
H^{0}(X_{1},{\cal O}_{X_{1}}(mK_{X})\otimes {\cal I}(h^{m}\mid_{X_{1}})
\cdot {\cal M}_{x_{1}}^{\lceil\sqrt[n_{1}]{\mu_{1}}(1-\varepsilon^{\prime})m\rceil})\neq 0
\]
holds.
$\square$ \end{lemma}

Let us  take  a general nonzero element $\sigma_{1,x_{1}}^{\prime}$ in
\[
H^{0}(X_{1},{\cal O}_{X_{1}}(m_{1}K_{X})\otimes{\cal I}(h^{m_{1}}\mid_{X_{1}})\cdot 
{\cal M}_{x_{1}}^{\lceil\sqrt[n_{1}]{\mu_{1}}(1-\varepsilon )m_{1}
\rceil}),
\]
for a sufficiently large $m_{1}$.
Using Lemma 3.6, let $\ell_{1}$ be as in Lemma 3.3 and let $\tau$ be a general nonzero
section in $H^{0}(X,{\cal O}_{X}(L_{1}))$, where $L_{1}$ is the line bundle 
as in Lemma 3.3. 
By Lemma 3.3, we may extend $\sigma_{1,x_{1}^{\prime}}\otimes\tau$ to 
a section 
\[
\sigma_{1,x_{1}}
\in H^{0}(X,{\cal O}_{X}((m_{1}+\ell_{1})K_{X})).
\]
As in Case 1.1, taking $\ell_{1}$ sufficiently large, we may assume that 
there exists a neighbourhood $U_{x_{1}}$ of $x_{1}$ such that 
$(\sigma_{1,x_{1}})$ is smooth on a  $U_{x_{1}}\backslash X_{1}$. 
We set 
\[
D_{1}(x_{1}) = \frac{1}{m_{1}+\ell_{1}}(\sigma_{1,x_{1}}). 
\]
Let $X_{1,\mbox{reg}}$ denote the regular locus of $X_{1}$. 
We may construct the divisors $\{ D_{1}(x_{1})\}$ as an algebraic 
family over $X_{1,\mbox{reg}}$.
Letting $x_{1}$ tend to $x$,
we obtain a {\bf Q}-divisor $D_{1}$ on $X$ which is $(m_{1}+\ell_{1})^{-1}$-times a divisor of a global holomorphic section
\[
\sigma_{1}\in H^{0}(X,{\cal O}_{X}((m_{1}+\ell_{1})K_{X})).
\]
By the construction, we may and do assume that there exists a neighbourhood $U_{x}$ of $x$ such that $(\sigma_{1})$ is smooth on $U_{x} \backslash X_{1}$. 
Let $\varepsilon_{0}$ be a sufficiently small 
positive rational number with $\varepsilon_{0} < \alpha_{0}$
such that  $(\alpha_{0}-\varepsilon_{0})D_{0}$ is not KLT at $x^{\prime}$
(this is possible because we are considering Case 2).

And we define the positive numbers $\alpha_{1}(x_{1})$ and $\alpha_{1}$ by
\[
\alpha_{1}(x_{1}) := \inf\{\alpha > 0 \mid 
\mbox{ $(\alpha_{0}-\varepsilon_{0})D_{0} + \alpha D_{1}(x_{1})$ 
is not KLT at $x_{1}$} \}.
\]
and
\[
\alpha_{1} := \inf\{\alpha > 0 \mid 
\mbox{ $(\alpha_{0}-\varepsilon_{0})D_{0} + \alpha D_{1}$ 
is  KLT at neither $x$ nor $x^{\prime}$} \}
\]
respectively. 
The definition of $\alpha_{1}$ is the same as in Case 1.1. 
But we note that $(\alpha_{0}-\varepsilon_{0})D_{0}$ is already not KLT at $x^{\prime}$.
We shall estimate $\alpha_{1}$. 
The proof of the following lemma is similar to that of Lemma 3.3. 
\begin{lemma}
Let $\delta$ be the fixed positive number as above.  Then we may assume that 
\[
\alpha_{1}\leqq \frac{n_{1}}{\sqrt[n_{1}]{\mu_{1}}} 
+ \delta 
\]
holds, if we take $\varepsilon^{\prime}$,  $\ell_{1}/m_{1}$ 
and $\varepsilon_{0}$ sufficiently 
small. $\square$
\end{lemma}
This estimate is better than Lemma 3.3. 
Then we may define the proper subvariety $X_{2}$ of $X_{1}$ 
as the minimal center of log canonical singularities 
of $(X,(\alpha_{0}-\varepsilon_{0})D_{0} + \alpha_{1} D_{1})$
at $x$ or $x^{\prime}$ as we have defined $X_{1}$. \vspace{5mm} \\

Lastly in Case 1.2 the construction of the 
filtration reduces to Case 2 as follows.
In Case 1.2, $X_{1}$ 
does not pass through $x^{\prime}$.  
Hence in this case the minimal center of LC singularities $X_{1}^{\prime}$ at 
$x^{\prime}$ does not pass through $x$. 
One may reduce Case 1.2 to Case 2, by ``strengthening'' the singularity of 
$D_{0}$ along $X_{1}^{\prime}$ as follows. 

Let $a_{1}$ be a sufficiently large positive integer such that 
\[
H^{0}(X,{\cal O}_{X}(a_{1}K_{X})\otimes {\cal I}_{X_{1}^{\prime}})
\neq 0.
\]
Let $\tau^{\prime}$ be a general nonzero section of
$H^{0}(X,{\cal O}_{X}(a_{1}K_{X})\otimes {\cal I}_{X_{1}^{\prime}})$. 
We note that there exists an effective {\bf Q}-divisor $G$ on $X$ such that 
$K_{X}-G$ is ample and $x$ is not contained in $\mbox{Supp}\, G$
as we have seen before.  Hence if we take $a_{1}$ sufficiently large, we may assume that the divisor
$(\tau^{\prime})$ does not contain $x$.  
In this case instead of $\sigma_{0}$, we shall use $\sigma_{0}^{e}\otimes \tau^{\prime}$, 
taking a positive integer $e$ large.  
Let $D_{0}^{\prime}:= (m_{0}e+a_{1})^{-1}(\sigma_{0}^{e}\otimes \tau^{\prime})$.Let us define a positive rational number $\alpha_{0}^{\prime}$ 
for $(X,D_{0}^{\prime})$ similar to $\alpha_{0}$. 
Then by the construction of $\tau^{\prime}$, then the minimal center of LC singularities of $(X,\alpha_{0}^{\prime}D_{0}^{\prime})$ at $x$ is 
$X_{1}$ and  $(X,\alpha_{0}^{\prime}D_{0}^{\prime})$ is not LC at $x^{\prime}$.
Also we can make $\alpha_{0}^{\prime}$ arbitrary close to $\alpha_{0}$ 
by taking $e$ sufficiently large.  
Hence we may assume that $\alpha_{0}^{\prime}$ satisfies the same 
estimate :
\[
\alpha_{0}^{\prime} \leqq \frac{n\sqrt[n]{2}}{\sqrt[n]{\mu_{0}}} + \delta
\]
as $\alpha_{0}$. 
And we may continue the construction of the filtration. 
In this way we can reduce Case 1.2 to Case 2.
\vspace{5mm} \\

In any case we construct the next stratum $X_{2}$ as the minimal center of log canonical singularities of $(X,(\alpha_{0}-\varepsilon_{0})D_{0} + \alpha_{1} D_{1})$
at $x$. 
If $X_{2}$ is a point, then we stop the construction of the filtration.
If $X_{2}$ is not a point, we continue exactly the same procedure 
replacing $X_{1}$ by $X_{2}$.  
And we continue the procedure as long as the new center of log canonical singularities ($X_{1},X_{2},\cdots$) is not a point. 
As a result, for any distinct points $x,x^{\prime}\in X^{\circ}$, inductively we construct 
a strictly decreasing sequence of subvarieties 
\[
X = X_{0}\supset X_{1}\supset \cdots \supset X_{r}\supset X_{r+1} =  x\,\,\mbox{or}\,\, x^{\prime}
\]
and invariants  :
\[
\alpha_{0} ,\alpha_{1},\ldots ,\alpha_{r},
\]
\[
\varepsilon_{0},\varepsilon_{1},\ldots ,\varepsilon_{r-1},
\]
\[
n >  n_{1}> \cdots > n_{r} \hspace{5mm}(n_{i} = \dim X_{i} ,i=1,\cdots ,r),
\]
and  
\[
\mu_{0},\mu_{1},\ldots ,\mu_{r}\hspace{5mm} (\mu_{i} := \mu (X_{i},(K_{X},h)\mid_{X_{i}}))
\]
depending on small positive  rational numbers $\varepsilon_{0},\ldots ,
\varepsilon_{r-1}$, large positive integers $m_{0},m_{1},\ldots ,m_{r}$,
positive integers $0=:\ell_{0}, \ell_{1},\cdots ,\ell_{r}$, 
\[
\sigma_{i}\in H^{0}(X,{\cal O}_{X}((m_{i}+\ell_{i})K_{X})) \hspace{10mm} (i= 0, \cdots ,r),
\]
\[
D_{i} = \frac{1}{m_{i}+ \ell_{i}}(\sigma_{i} ) \hspace{10mm}(i= 0, \cdots ,r),
\]
etc.

By Nadel's vanishing theorem (\cite[p.561]{n}) we have the following lemma.

\begin{lemma} 
For every positive integer  $m > 1+ \sum_{i=0}^{r}\alpha_{i}$, 
$\Phi_{\mid mK_{X}\mid}$ separates $x$ and $x^{\prime}$.
And we may assume that 
\[
\alpha_{i}\leqq \frac{n_{i}\sqrt[n_{i}]{2}}{\sqrt[n_{i}]{\mu_{i}}} + \delta 
\]
holds for  every $0\leqq i\leqq r$.
$\square$ \end{lemma}
{\bf Proof}.
For $i= 0,1,\ldots, r$, let $h_{i}$ be the singular hermitian metric 
on $K_{X}$ defined by 
\[
h_{i}:= \frac{1}{\mid\sigma_{i}\mid^{\frac{2}{m_{i}+\ell_{i}}}},
\]
where we have set $\ell_{0}= 0$.   
Using Kodaira's lemma (\cite[Appendix]{k-o}), 
let us take an effective {\bf Q}-divisor $G$ on $X$ such that
$K_{X} - G$ is ample as before.  As before we may assume that 
$\mbox{Supp}\, G$ contains neither $x$ nor $x^{\prime}$.
Let $h^{\prime}_{G}$ be a $C^{\infty}$-hermitian metric on the {\bf Q}-line bundle
$K_{X}- G$ with strictly positive curvature. 
Let $G = \sum_{k}g_{k}G_{k}$ be the irreducible decomposition of $G$ 
and let $\sigma_{G_{k}}$ be a global holomorphic section of 
${\cal O}_{X}(G_{k})$ with divisor $G_{k}$. 
Then  
\[
h_{G}:= h_{G}^{\prime}\cdot 
(\prod_{k}\frac{1}{\mid\sigma_{G_{k}}\mid^{2g_{k}}})
\]
is a singular hermitian metric of $K_{X}$ with strictly positive curvature
current.

Let $m$ be a positive integer such that $m > 1+ \sum_{i=0}^{r}\alpha_{i}$
as above. 
Let $\varepsilon_{G}$ be a positive number such that 
\[
\varepsilon_{G} < m - 1 - ( \sum_{i=0}^{r-1}(\alpha_{i}-\varepsilon_{i}) +\alpha_{r}).
\] 
We set 
\[
\beta := \sum_{i=0}^{r-1}(\alpha_{i}-\varepsilon_{i}) +\alpha_{r} + \varepsilon_{G}. 
\] 
\[
h_{x,x^{\prime}} = (\prod_{i=0}^{r-1}h_{i}^{\alpha_{i}-\varepsilon_{i}})\cdot 
 h_{r}^{\alpha_{r}}\cdot h^{m-1-\beta}\cdot h_{G}^{\varepsilon_{G}}.
\]
Then we see that  ${\cal I}(h_{x,x^{\prime}})$ defines a subscheme of 
$X$ with isolated support around $x$ or $x^{\prime}$ by the definition of 
the invariants $\{\alpha_{i}\}$'s and the fact that $\mbox{Supp}\, G$ 
contains neither  $x$ nor $x^{\prime}$. 
By the construction the curvature current $\Theta_{h_{x,x^{\prime}}}$ is strictly positive on $X$. 
Then by Nadel's vanishing theorem (\cite[p.561]{n}) we see that 
\[
H^{1}(X,{\cal O}_{X}(mK_{X})\otimes {\cal I}(h_{x,x^{\prime}})) = 0
\]
holds. 
Hence 
\[
H^{0}(X,{\cal O}_{X}(mK_{X}))
\rightarrow 
H^{0}(X,{\cal O}_{X}(mK_{X})\otimes {\cal O}_{X}/{\cal I}(h_{x,x^{\prime}}))
\]
is surjective. 
Since by the construction of $h_{x,x^{\prime}}$,  
$\mbox{Supp}({\cal O}_{X}/{\cal I}(h_{x,x^{\prime}}))$ 
contains both $x$ and $x^{\prime}$ and is 
isolated at least at  one of $x$ or $x^{\prime}$. 
Hence by the above surjection, there exists a section
$\sigma\in H^{0}(X,{\cal O}_{X}(mK_{X}))$ such that 
\[
\sigma (x) \neq 0,\sigma (x^{\prime}) = 0
\]
or 
\[
\sigma (x) = 0,\sigma (x^{\prime}) \neq 0
\]
holds. 
This implies that $\Phi_{\mid mK_{X}\mid}$ separates 
$x$ and $x^{\prime}$.  
The proof of the last statement is similar to  that of Lemma 3.3.  
$\square$.

\subsection{Estimate of the degree}

To relate $\mu_{0}$ and the degree of pluricanonical images of $X$, 
we need the following lemma. 

\begin{lemma}
If $\Phi_{\mid mK_{X}\mid}$ is a birational rational map
onto its image, then
\[
\deg \Phi_{\mid mK_{X}\mid}(X)\leqq \mu_{0}\cdot m^{n}
\]
holds.
$\square$ \end{lemma}
{\bf Proof}.
Let $p : \tilde{X}\longrightarrow X$ be the resolution of 
the base locus of $\mid mK_{X}\mid$ and let 
\[
p^{*}\mid mK_{X}\mid = \mid P_{m}\mid + F_{m}
\]
be the decomposition into the free part $\mid P_{m}\mid$ 
and the fixed component $F_{m}$. 
We have
\[
\deg \Phi_{\mid mK_{X}\mid}(X) = P_{m}^{n}
\]
holds.
Then by the ring structure of $R(X,K_{X})$, we have an injection 
\[
H^{0}(\tilde{X},{\cal O}_{\tilde{X}}(\nu P_{m}))\rightarrow 
H^{0}(X,{\cal O}_{X}(m\nu K_{X})\otimes{\cal I}(h^{m\nu}))
\]
for every $\nu\geqq 1$, since 
the righthand side is isomorphic to 
$H^{0}(X,{\cal O}_{X}(m\nu K_{X}))$ by the definition of 
an AZD.
We note that since ${\cal O}_{\tilde{X}}(\nu P_{m})$ is globally generated
on $\tilde{X}$, for every $\nu \geqq 1$ we have the injection 
\[
{\cal O}_{\tilde{X}}(\nu P_{m})\rightarrow p^{*}({\cal O}_{X}(m\nu K_{X})\otimes{\cal I}(h^{m\nu})).
\]
Hence there exists a natural homomorphism 
\[
H^{0}(\tilde{X},{\cal O}_{\tilde{X}}(\nu P_{m}))
\rightarrow 
H^{0}(X,{\cal O}_{X}(m\nu K_{X})\otimes{\cal I}(h^{m\nu}))
\]
for every $\nu\geqq 1$. 
This homomorphism is clearly injective. 
This implies that 
\[
\mu_{0} \geqq  m^{-n}\cdot\mu (\tilde{X},P_{m})
\]
holds by the definition of $\mu_{0}$. 
Since $P_{m}$ is nef and big on $X$, we see that 
\[
\mu (\tilde{X},P_{m}) = P_{m}^{n}
\]
holds.
Hence
\[
\mu_{0}\geqq m^{-n}\cdot P_{m}^{n}
\]
holds.  This implies that
\[
\deg \Phi_{\mid mK_{X}\mid}(X)\leqq \mu_{0}\cdot m^{n}
\]
holds.
$\square$

\subsection{Use of the subadjunction theorem}

Let 
\[
X = X_{0}\supset X_{1}\supset \cdots \supset X_{r}\supset X_{r+1} =  x\,\,\mbox{or}\,\, x^{\prime}
\]
be the filtration of $X$ as in Section 3.1.
\begin{lemma}
Let $W_{j}$ be  a nonsingular model of $X_{j}$.
For every  $W_{j}$, 
\[
\mu (W_{j},K_{W_{j}})
\leqq 
(\lceil (1+\sum_{i=0}^{j-1}\alpha_{i})\rceil )^{n_{j}}\cdot\mu_{j}
\]
holds, where $\mu_{j} = (K_{X},h)^{n_{j}}\cdot X_{j}$ as in Section 3.1 (
we note that $\mu (W_{j},K_{W_{j}})$ depends only on $X_{j}$).
$\square$ \end{lemma}
{\bf Proof}.

Let us set 
\[
\beta_{j}:= \varepsilon_{j-1}+\sum_{i=0}^{j-1} (\alpha_{i}-\varepsilon_{i}).
\]
Let $D_{i}$ denote the divisor $m_{i}^{-1}(\sigma_{i})$ and we set 
\[
D : = \sum_{i=1}^{j-1} (\alpha_{i}-\varepsilon_{i})D_{i}+\varepsilon_{j-1}D_{j-1}.
\]
Let $\pi : Y \longrightarrow X$ be a log resolution of 
$(X,D)$ 
which factors through an embedded resolution 
$\varpi : W_{j} \longrightarrow X_{j}$ of $X_{j}$. 
By the modification as in Section 3.1, we may assume that there exists a unique  irreducible component $F_{j}$ of the exceptional divisor with discrepancy $-1$ which dominates $X_{j}$. 
Let 
\[
\pi_{j} : F_{j} \longrightarrow W_{j}
\]
be the natural morphism induced by the construction. 
We set 
\[
\pi^{*}(K_{X} + D)\mid_{F_{j}} = K_{F_{j}} + G.
\]
We may assume that the support of $G$ is a divisor with normal crossings.
Then all the coefficients of the horizontal component $G^{h}$ of $G$ with respect to $\pi_{j}$ are less than $1$  because $F_{j}$ is the unique exceptional divisor with discrepancy $-1$. 

Let $dV$ be a $C^{\infty}$-volume form on the  $X$. 
Let $\Psi$ be the function defined by 
\[
\Psi := \log (h^{\beta_{j}}\cdot\mid\sigma_{j-1}\mid^{{\frac{2\varepsilon_{j-1}}{m_{j-1}}}}\cdot\prod_{i=0}^{j-1}\mid \sigma_{i}\mid^{\frac{2(\alpha_{i}-\varepsilon_{i})}{m_{i}}}).
\]
Then the residue $\mbox{Res}_{F_{j}}(\pi^{*}(e^{-\Psi}\cdot dV))$ of  $\pi^{*}(e^{-\Psi}\cdot dV)$
to $F_{j}$ 
is a singular volume form with algebraic singularities corresponding 
to the divisor $G$. 
Since every coefficient of $G^{h}$ is less than $1$, there exists a nonempty Zariski open subset 
$W_{j}^{0}$ of $W_{j}$ such that 
 $\mbox{Res}_{F_{j}}(\pi^{*}(e^{-\Psi}\cdot dV))$ is integrable
 on $\pi_{j}^{-1}(W_{j}^{0})$.

Then  the pullback of the residue $dV[\Psi ]$ of $e^{-\Psi}\cdot dV$ 
(to $X_{j}$) to 
 $W_{j}$ is given by the fiber integral of 
the above singular volume form $\mbox{Res}_{F_{j}}(\pi^{*}(e^{-\Psi}\cdot dV))$ on $F_{j}$, i.e.,
\[
\varpi^{*}dV[\Psi ] = \int_{F_{j}/W_{j}}\mbox{Res}_{F_{j}}(\pi^{*}(e^{-\Psi}\cdot dV))
\]
holds. 
By Theorem \ref{pos}, we see that 
$(K_{F_{j}}+G) -\pi_{j}^{*}(K_{W_{j}}+\Delta )$ is nef,
where $\Delta$ is the {\bf Q}-divisor defined as in Theorem \ref{pos}. 
We note that $K_{F_{j}}+G$ is {\bf Q}-linear equivalent to 
$(1+\beta_{j})\pi^{*}K_{X}$ by the construction. 
Hence we see that $(1+\beta_{j})\varpi^{*}K_{X} - (K_{W_{j}}+\Delta )$ 
is nef and 
\begin{equation}\label{delta}
\mu (W_{j},K_{W_{j}})\leqq \mu (W_{j},(1+\beta_{j})\varpi^{*}(K_{X}\mid_{X_{j}})-\Delta )
\end{equation}
holds.

Let $e$ be a positive integer such that $e\cdot \Delta$ is an integral divisor.
Let $\sigma_{e\cdot\Delta}$ be a meromorphic section of ${\cal O}_{W_{j}}(e\cdot\Delta )$
with divisor $e\cdot\Delta$.
Then we may consider the $e$-th root $\sigma_{\Delta}$ of $\sigma_{e\cdot\Delta}$ 
as a multivalued  meromorphic section of the 
{\bf Q}-line bundle ${\cal O}_{W_{j}}(\Delta )$ with divisor $\Delta$.
Let $h_{\Delta}$ be a $C^{\infty}$-hermitian metric on the {\bf Q}-line bundle
${\cal O}_{W_{j}}(\Delta )$, i.e.,  $h_{\Delta}$ is the  $e$-th root of 
a $C^{\infty}$-hermitian metric on the line bundle ${\cal O}_{W_{j}}(e\cdot\Delta )$.  Then $h_{\Delta}(\sigma_{\Delta},\sigma_{\Delta})$ is a single valued 
funtion on $W_{j}$. 

Let us recall the interpretation of the divisor $\Delta$ in 
Section 3.7.
Let $dV_{W_{j}}$ be a $C^{\infty}$-volume form on $W_{j}$.
We note that in the above definition of the function $\Psi$, we have used 
$h^{\beta_{j}}$ instead of $dV^{-\beta_{j}}$.  
Hence we see that there exists a positive constant $C$ such that 
\begin{equation}\label{fiber}
\varpi^{*}dV[\Psi ] = \int_{F_{j}/W_{j}}\mbox{Res}_{F_{j}}(\pi^{*}(e^{-\Psi}\cdot dV))
\leqq C\cdot\frac{\varpi^{*}(dV\cdot h)^{-\beta_{j}}}{h_{\Delta}(\sigma_{\Delta},\sigma_{\Delta})}\cdot dV_{W_{j}}
\end{equation}
hold.

We may assume that $\beta_{j}$ is not an integer without loss of generality. 
In fact this can be satisfied, if we perturb $\varepsilon_{0},\cdots, 
\varepsilon_{j-1}$ or $\sigma_{0},\cdots ,\sigma_{j-1}$. 
And passing to the lmit, the general case follows. 
This condition is  to assure the inequality  
$\lceil 1 + \beta_{j}\rceil > 1 +\beta_{j}$ 
and  this inequality corresponds to the condition : $d > \alpha m_{0}$ in Theorem \ref{subad1}.
We note that for every positive integer $m$, 
every global holomorphic section of $mK_{X}$ is 
bounded with respect to $h^{m}$.  
Then since  the curvature current $\Theta_{h}$ is semipositive in the sense of current, applying  Theorem \ref{subad1} (see also Remark \ref{r2.6}
for the selfcontainedness) , we have the interpolation :
\[
A^{2}(W_{j},m(\lceil 1+\beta_{j}\rceil )\varpi^{*}K_{X}, 
\varpi^{*}(e^{-(m-1)\varphi}\cdot dV^{-m}\otimes h^{m\lceil \beta_{j}\rceil}),
\varpi^{*}dV[\Psi ])
\]
\[
\hspace{60mm} 
\rightarrow 
H^{0}(X,{\cal O}_{X}(m(\lceil 1+ \beta_{j}\rceil )K_{X})),
\]
where $\varphi$ is the weight function defined by 
\[
\varphi := \log \frac{dV_{W_{j}}}{\varpi^{*}dV[\Psi]}
\]
as in Theorem \ref{subad1}. 
By (\ref{fiber}), we see that
\begin{equation}\label{varphi}
\varpi^{*}(e^{-\varphi}\cdot dV^{-1}\otimes h^{\beta_{j}}\mid_{X_{j}})
\leqq C\cdot h_{\Delta}(\sigma_{\Delta},\sigma_{\Delta})^{-1}\cdot \varpi^{*}(dV^{-(1+\beta_{j})}\mid_{X_{j}})
\end{equation}
holds.  We note that $\Delta$ may not be effective. Hence a priori 
the element of $A^{2}(W_{j},m(\lceil 1+\beta_{j}\rceil )\varpi^{*}K_{X}, 
\varpi^{*}(e^{-(m-1)\varphi}\cdot dV^{-m}\otimes h^{m\lceil \beta_{j}\rceil}),
\varpi^{*}dV[\Psi ])$ may have pole along the degenerate locus (zero locus) of 
$\varpi^{*}dV[\Psi ]$. 
But this cannot occur by the existence of the extension 
and the birational invariance of plurigenera. 
As in the remark below, we also may reduce the proof to the case that 
$\Delta$ is effective. 

Since   
$(1+\beta_{j})\varpi^{*}(K_{X}\mid_{X_{j}})-(K_{W_{j}}+\Delta )$ is 
nef (This is  nothing but the main part  of the proof of 
Kawamata's subadjunction theorem 
\cite[Theorem 1]{ka}.  Then the proof of \cite[Theorem 1]{ka} follows from the 
observation that $\varpi_{*}\Delta$ is effective), 
by using Theorem \ref{pos}, noting the equality $dV[\Psi ] = e^{-\varphi}\cdot dV_{W_{j}}$, 
the inequalities (\ref{delta}),(\ref{varphi}) and the existence of the above interpolation
imply that 
\begin{eqnarray*}
\mu(W_{j},K_{W_{j}}) \leqq \hspace{120mm} \\
 n_{j}!\cdot\overline{\lim}_{m\rightarrow\infty}
m^{-n_{j}}\dim \mbox{Image}\, \{H^{0}(X,{\cal O}_{X}(m(\lceil 1+ \beta_{j}\rceil )K_{X}))
\rightarrow H^{0}(X_{j},{\cal O}_{X_{j}}(m(\lceil 1+ \beta_{j}\rceil )K_{X}))\} 
\end{eqnarray*}
holds.  Here we have used the fact that for any pseudoeffective divisors 
$M_{1},M_{2}$ on a smooth projective variety $V$ such that $M_{1}- M_{2}$ is pseudoeffective, the inequality: $\mu (V,M_{1})\geqq \mu (V,M_{2})$ holds
(the proof is trivial and left to the reader). 

Since  every element of 
$H^{0}(X,{\cal O}_{X}(m(\lceil 1+ \beta_{j}\rceil )K_{X}))$ 
is  bounded on $X$  with respect to $h^{m(\lceil 1+\beta_{j}\rceil )}$
(cf. Remark \ref{r2.3}). 
In particular the restriction of an element of $H^{0}(X,{\cal O}_{X}(m(\lceil 1+ \beta_{j}\rceil )K_{X}))$ to $X_{j}$ is bounded with respect to 
$h^{m(\lceil 1+\beta_{j}\rceil )}\mid_{X_{j}}$.
Hence by the existence of the above interpolation, we have that 
\begin{equation}
\mu (W_{j},K_{W_{j}})
\leqq 
\mu (X_{j},(\lceil 1+\beta_{j}\rceil )K_{X},
h^{\lceil 1+\beta_{j}\rceil})\mid_{X_{j}})
\end{equation}
holds. 
This is the only point where Theorem \ref{subad1} is used. 

By the trivial inequality 
\[
\beta_{j} \leqq \sum_{i=0}^{j-1}\alpha_{i}.
\]
we have that 
\[
\mu (W_{j},K_{W_{j}})\leqq (\lceil 1+\sum_{i=0}^{j-1}\alpha_{i}\rceil )^{n_{j}}(K_{X},h)^{n_{j}}\cdot X_{j}\]
holds by the definition of $(K_{X},h)^{n_{j}}\cdot X_{j}$.
This is the desired inequality, since $\mu_{j} = (K_{X},h)^{n_{j}}\cdot X_{j}$ 
holds by the definition of $\mu_{j}$.
$\square$ 
\begin{remark}
In the above proof, the divisor $\Delta$ on $W_{j}$ may not be effective. 
But it is clear that $\varpi_{*}\Delta$ is effective (cf. the proof of \cite[Theorem 1]{ka}).  
If we replace $X_{j}$ by $W_{j}$ and $X$ by the the ambient space of the 
embedded resolution $\varpi : W_{j}\longrightarrow X_{j}$, 
we may reduce the above proof to the case that $X_{j}$ is already smooth. 
In this case we may assume that $\Delta$ is  effective. 
\end{remark}

Now we shall complete the proofs of Theorems 1.1 and 1.2.

Suppose that Theorem 1.2 holds for every projective varieties of general 
type of dimension $< n$, i.e., there exist positive constants $\{C(k) (k <n)\}$
such that for  every smooth projective $k$-fold $Y$ of general type
\[
\mu (Y,K_{Y}) \geqq C(k)
\] 
holds. 
Let $X$ be a smooth projective variety of general type of 
dimension $n$. 
Let $U_{0}$ be a nonempty open subset of $X$ with respect to 
{\bf countable Zariski topology} such that for every $x\in U_{0}$ 
there exist no subvarieties of nongeneral type containing $x$. 
Such a set $U_{0}$ surely exists, since there exists no dominant family of 
 subvarieties of nongeneral type in $X$.
 In fact if such a dominant family exists, then this contradicts the
 assumption that $X$ is of general type.   
Then if $(x,x^{\prime})\in (U_{0}\times U_{0})\backslash\Delta_{X}$, the stratum $X_{j}$ 
as in Section 3.1 is of general type for every $j$ by the definition of $U_{0}$. 
By Lemma 3.10 and the definition of $C(n_{j})$, we see that
\begin{equation}
C(n_{j}) \leqq (\lceil (1+\sum_{i=0}^{j-1}\alpha_{i})\rceil )^{n_{j}}\cdot \mu_{j}
\end{equation}
holds for $W_{j}$. 
Since 
\begin{equation}
\alpha_{i} \leqq \frac{\sqrt[n_{i}]{2}\, n_{i}}{\sqrt[n_{i}]{\mu_{i}}}+ \delta
\end{equation}
holds for every $0\leqq i\leqq r$ by Lemma 3.8, combining (5) and (6), 
we see that  
\[
\frac{1}{\sqrt[n_{j}]{\mu_{j}}}\leqq 
(2+\sum_{i=0}^{j-1}\frac{\sqrt[n_{i}]{2}\, n_{i}}{\sqrt[n_{i}]{\mu_{i}}} )\cdot C(n_{j})^{-\frac{1}{n_{j}}}
\]
holds for every $j \geqq 1$.

Using the above inequality inductively, we obtain the following lemma. 
\begin{lemma}
Suppose that $\mu_{0} \leqq 1$ holds.
Then there exists a positive constant $C$ depending only on $n$ 
such that for every $(x,x^{\prime})\in (U_{0}\times U_{0})\backslash\Delta_{X}$ 
the corresponding invariants $\{ \mu_{0},\cdots ,\mu_{r}\}$ 
and $\{ n_{1},\cdots ,n_{r}\}$ depending 
on $(x,x^{\prime})$ ($r$ may also depend on $(x,x^{\prime})$) satisfies
the inequality :
\[
2+\lceil \sum_{i=0}^{r}\frac{\sqrt[n_{i}]{2}\, n_{i}}{\sqrt[n_{i}]{\mu_{i}}}\rceil \leqq \lfloor\frac{C}{\sqrt[n]{\mu_{0}}}\rfloor .
\]
$\square$ \end{lemma} 
We note that  $\{ n_{1},\cdots ,n_{r}\}$ is a 
strictly decreasing sequence and this sequence has only finitely many 
possibilities. 
By Lemmas 3.8 and 3.11  we see that 
for 
\[
m := \lfloor\frac{C}{\sqrt[n]{\mu_{0}}}\rfloor ,
\]
$\mid mK_{X}\mid$ separates points on $U_{0}$.
Hence  $\mid mK_{X}\mid$ gives a birational embedding of $X$.

Then by Lemma 3.9, if $\mu_{0} \leqq 1$ holds, 
\[
\deg \Phi_{\mid mK_{X}\mid}(X) 
\leqq C^{n}
\]
holds. 
Also 
\[
\dim H^{0}(X,{\cal O}_{X}(mK_{X}))
\leqq n+1 + \deg \Phi_{\mid mK_{X}\mid}(X) 
\]
holds by the semipositivity of the $\Delta$-genus (\cite{fu}). 
Hence we have that if $\mu_{0} \leqq 1$,
\[
\dim H^{0}(X,{\cal O}_{X}(mK_{X})) 
\leqq n+1+ C^{n} 
\]
holds. 

Since $C$ is a positive constant depending only on $n$, 
combining  the above two inequalities, 
we have that there exists a positive constant $C(n)$ 
depending only on $n$ such that 
\[
\mu_{0}  \geqq C(n)
\]
holds. 

More precisely we argue as follows. 
Let ${\cal H}$ be the union of the irreducible components of the Hilbert scheme of  projective spaces of dimension $\leqq n+ C^{n}$
and the degree $\leqq C^{n}$.
By the general theory of Hilbert schemes (\cite[expos\'{e} 221]{gro}), ${\cal H}$ consists of finitely many irreducible components.  
Let ${\cal H}_{0}$ be the Zariski open subset of ${\cal H}$ which parametrizes 
irreducible subvarieties. 
Then there exists a finite stratification of ${\cal H}_{0}$ by Zariski locally closed subsets such that on each stratum, there exists a simultaneous resolution 
of the universal family on the stratum. 
We note that the volume of the canonical bundle of the resolution is constant on each stratum by \cite{tu6,nak}.
Hence there exists a positive constant $C(n)$ depending only on $n$  such that 
\[
\mu (X,K_{X})\geqq C(n)
\]
holds for every projective $n$-fold $X$ of general type 
by the degree bound as above. 
This completes the proof of Theorem 1.2. \vspace{5mm}$\square$\\

Now let us prove Theorem 1.1.
Then by Lemmas 3.8 and 3.11, we see that there exists 
a positive integer $\nu_{n}$ depending only on $n$ such that 
for every projective $n$-fold $X$ of general type, 
$\mid mK_{X}\mid$ gives a birational embedding into a 
projective space for every $m\geqq \nu_{n}$. 
This completes the proof of Theorem 1.1. $\square$

\section{The Severi-Iitaka conjecture}

Let $X$ be a smooth projective variety.
We set
\[
Sev (X) := \{ (f,[Y])\mid f : X\longrightarrow Y\,\,\,\,\mbox{dominant
rational map and $Y$ is of general type}\} ,
\]
where $[Y]$ denotes the birational class of $Y$.
By Theorem 1.1 and \cite[p.117, Proposition 6.5]{m}
we obtain the following theorem.

\begin{theorem}
$Sev (X)$ is finite.
$\square$ \end{theorem}
\begin{remark}
In the case of $\dim Y = 1$, Theorem 4.1 is known 
as Severi's theorem. 
In the case of $\dim Y = 2$, Theorem 4.1 has already been known 
by K. Maehara (\cite{m}).
In the case of $\dim Y = 3$, Theorem 4.1 has recently proved by
T. Bandman and G. Dethloff (\cite{b-d}).
$\square$ \end{remark}

Author's address\\
Hajime Tsuji\\
Department of Mathematics\\
Sophia University\\
7-1 Kioicho, Chiyoda-ku 102-8554\\
Japan \\
e-mail address: tsuji@mm.sophia.ac.jp

\end{document}